\documentclass[12pt,a4paper]{article}
\usepackage{etoolbox} %
\usepackage{amsfonts,amssymb}
\usepackage[margin=1in]{geometry}
\usepackage{microtype}
\usepackage{times}
\usepackage[amsthm,thmmarks]{ntheorem}
\usepackage{amsmath,xypic}
\usepackage{multicol}
\usepackage{hyphenat}
\usepackage[usenames,dvipsnames]{color}
\usepackage{rotating}
\usepackage{longtable}
\usepackage{caption}
\usepackage{tcolorbox}
\usepackage{stackengine}

\usepackage{epigraph}
\usepackage{enumitem}
\setlist[enumerate]{listparindent=0.5in}
\usepackage{mathrsfs}
\DeclareMathAlphabet{\mathscrbf}{OMS}{mdugm}{b}{n} %
\newcommand{\be}{\begin{equation}}
\newcommand{\ee}{\end{equation}}
\newcommand{\bes}{\begin{equation*}}
\newcommand{\ees}{\end{equation*}}
\newcommand{\bea}{\begin{eqnarray}}
\newcommand{\eea}{\end{eqnarray}}
\newcommand{\beas}{\begin{eqnarray}}
\newcommand{\eeas}{\end{eqnarray}}
\newcommand{\ben}{\begin{note}}
\newcommand{\een}{\end{note}}
\newcommand{\bexl}{\vskip0.1em\noindent\hrulefill\vskip1em\begin{ExerciseList}}
\newcommand{\eexl}{\end{ExerciseList}\hrulefill}

\newcommand{\bthm}{\begin{theorem}}
\newcommand{\ethm}{\end{theorem}}
\newcommand{\bpro}{\begin{prop}}
\newcommand{\epro}{\end{prop}}
\newcommand{\bcor}{\begin{corollary}}
\newcommand{\ecor}{\end{corollary}}
\newcommand{\bcon}{\begin{conjecture}}
\newcommand{\econ}{\end{conjecture}}
\newcommand{\bp}{\begin{proof}}
\newcommand{\ep}{\end{proof}}
\newcommand{\blem}{\begin{lemma}}
\newcommand{\elem}{\end{lemma}}
\newcommand{\bn}{\begin{note}}
\newcommand{\en}{\end{note}}
\newcommand{\benum}{\begin{enumerate}}
\newcommand{\eenum}{\end{enumerate}}
\newcommand{\bed}{\begin{defn}}
\newcommand{\eed}{\end{defn}}
\newcommand{\brem}{\begin{remark}}
\newcommand{\erem}{\end{remark}}

\newcommand{\btik}{\begin{tikzpicture}\begin{axis}[scale=0.5,axis y line=center, axis x line=middle]}
\newcommand{\etik}{\end{axis}\end{tikzpicture}}

\let\into=\hookrightarrow
\let\mapsto=\longmapsto

\newcommand{\upperRomannumeral}[1]{\uppercase\expandafter{\romannumeral#1}}

 \usepackage{tikz}
\usetikzlibrary{mindmap,backgrounds}
\usepackage{graphicx}
\usepackage[pagewise]{lineno}
\usepackage{stmaryrd}
\newtoggle{arxiv}
\toggletrue{arxiv}
\iftoggle{arxiv}{
	\usepackage{url}
	\usepackage{natbib}
	\let\cite=\citep
    \bibliographystyle{plainnat}
	\typeout{In the arxiv mode}
	\usepackage{fancyhdr}
	\usepackage[colorlinks,citecolor=blue]{hyperref}
}
{
	\usepackage[backend=biber,style=alphabetic,citestyle=authoryear]{biblatex} %
	\usepackage{fancyhdr}
	\usepackage[colorlinks,citecolor=blue,colorlinks=true,hyperindex, citecolor=blue, urlcolor=blue]{hyperref}
	\let\cite=\parencite
\addbibresource{../../master/masterofallbibs.bib}
}
\newtoggle{draft}
\toggletrue{draft}
\newtheorem{theorem}[equation]{Theorem}      %
\newtheorem{lemma}[equation]{Lemma}          %
\newtheorem{corollary}[equation]{Corollary}  %
\newtheorem{proposition}[equation]{Proposition}

\theoremstyle{definition}
\newtheorem{conj}[equation]{Conjecture}

\theoremstyle{definition}
\newtheorem{defn}[equation]{Definition}
\theoremstyle{remark}

\theoremsymbol{\ensuremath{\bullet}}

\theoremstyle{definition}
\newtheorem{remark}[equation]{Remark}
\theoremsymbol{\ensuremath{\bullet}}

\numberwithin{equation}{subsection}
\newcommand{\parat}[1]{\subsubsection{#1}}
\newcommand{\para}{\parat{}}

\let\into=\hookrightarrow
\let\isom=\simeq

\newcommand{\A}{\mathcal{A}}
\newcommand{\abs}[1]{\left\vert#1\right\vert}

\newcommand{\bF}{{\bar{F}}}
\newcommand{\bQ}{{\bar{\Q}}}

\newcommand{\C}{{\mathbb C}}

\newcommand{\F}{{\mathbb F}}

\newcommand{\gal}{{\rm Gal}}

\newcommand{\N}{\mathcal{N}}

\newcommand{\Q}{{\mathbb Q}}
\newcommand{\R}{{\mathbb R}}

\newcommand{\spec}{{\rm Spec}}
\newcommand{\Spec}{{\rm Spec}}

\newcommand{\Z}{{\mathbb Z}}

\renewcommand{\int}{\operatorname{int}}
\renewcommand{\O}{{\mathcal O}}
\renewcommand{\P}{{\mathbb P}}

\newcommand{\invlim}{\varprojlim}

\renewcommand{\bpro}{\begin{proposition}}
	\renewcommand{\epro}{\end{proposition}}
\renewcommand{\bcon}{\begin{conj}}
	\renewcommand{\econ}{\end{conj}}
\newcommand{\ilim}{\varprojlim}
\newcommand{\ENDDOCUMENT}{
\iftoggle{arxiv}{%
	\bibliography{../../master/masterofallbibs.bib}
}
{\printbibliography} %
\end{document}
}

\let\mathcal=\mathscr

\title{Untilts of fundamental groups: construction of labeled isomorphs of fundamental groups\\
{\large Arithmetic Holomorphic Structures}\\
\textcolor{red}{\normalsize{Preliminary Version For Comments}}}
\author{Kirti Joshi}

\begin{document}
\maketitle

\lhead{}

\epigraphwidth0.55\textwidth
\epigraph{Heard melodies are sweet, but those unheard \\ Are sweeter; therefore, ye soft pipes, play on;}{\iftoggle{arxiv}{John Keats \cite{keats}}{\citeauthor{keats}}}
\newcommand{\act}{\curvearrowright}
\newcommand{\lmp}{{\Pi\act\Ot}}
\newcommand{\lmpi}{{\lmp}_{\int}}
\newcommand{\lmpf}{\lmp_F}
\newcommand{\Om}{\O^{\times\mu}}
\newcommand{\Omf}{\O^{\times\mu}_{\bF}}
\renewcommand{\N}{\mathbb{N}}
\newcommand{\yoga}{Yoga}
\newcommand{\gl}[1]{{\rm GL}(#1)}
\newcommand{\bK}{\overline{K}}
\newcommand{\reptrip}{\rho:G_K\to\gl V}
\newcommand{\reptripp}[1]{\rho\circ\alpha:G_{\ifstrempty{#1}{K}{{#1}}}\to\gl V}
\newcommand{\benumlab}{\begin{enumerate}[label={{\bf(\arabic{*})}}]}
\newcommand{\ord}{\mathop{\rm ord}\nolimits}	
\newcommand{\kcs}{K^\circledast}
\newcommand{\lcs}{L^\circledast}
\renewcommand{\A}{\mathbb{A}}
\newcommand{\bfq}{\bar{\mathbb{F}}_q}
\newcommand{\tripod}{\P^1-\{0,1728,\infty\}}

\newcommand{\vseq}[2]{{#1}_1,\ldots,{#1}_{#2}}
\newcommand{\anab}[4]{\left({#1},\{#3 \}\right)\anabelmap\left({#2},\{#4 \}\right)}

\newcommand{\cpt}{\C_p^\flat}

\newcommand{\gln}{{\rm GL}_n}
\newcommand{\glo}[1]{{\rm GL}_1(#1)}
\newcommand{\glt}[1]{{\rm GL_2}(#1)}

\newcommand{\iut}{\cite{mochizuki-iut1, mochizuki-iut2, mochizuki-iut3,mochizuki-iut4}}
\newcommand{\topics}{\cite{mochizuki-topics1,mochizuki-topics2,mochizuki-topics3}}

\newcommand{\linv}{\mathfrak{L}}
\newcommand{\bedef}{\begin{defn}}
\newcommand{\eedef}{\end{defn}}
\renewcommand{\act}[1][]{\overset{#1}{\curvearrowright}}
\newcommand{\bfx}{\overline{F(X)}}
\newcommand{\anabelmap}{\leftrightsquigarrow}
\newcommand{\ban}[1][G]{\mathcal{B}({#1})}
\newcommand{\pit}{\Pi^{temp}}
 
 \newcommand{\bL}{\overline{L}}
 \newcommand{\bkm}{\bK_M}
 \newcommand{\vbk}{v_{\bK}}
 \newcommand{\vbkm}{v_{\bkm}}
\newcommand{\ocs}{\O^\circledast}
\newcommand{\ot}{\O^\triangleright}
\newcommand{\ocsk}{\ocs_K}
\newcommand{\otk}{\ot_K}
\newcommand{\ok}{\O_K}
\newcommand{\oko}{\O_K^1}
\newcommand{\oks}{\ok^*}
\newcommand{\Qpb}{\overline{\Q}_p}
\newcommand{\Qpbh}{\widehat{\overline{\Q}}_p}
\newcommand{\tr}{\triangleright}
\newcommand{\ocpt}{\O_{\C_p}^\tr}
\newcommand{\ocpf}{\O_{\C_p}^\flat}
\newcommand{\sG}{\mathscr{G}}
\newcommand{\sxqp}{\mathscr{X}_{\cpt,\Q_p}}
\newcommand{\syqp}{\mathscr{Y}_{\cpt,\Q_p}}
\newcommand{\sxfe}{\mathscr{X}_{F,E}}
\newcommand{\sxfep}{\mathscr{X}_{F,E'}}
\newcommand{\syfe}{\mathscr{Y}_{F,E}}
\newcommand{\syfep}{\mathscr{Y}_{F,E'}}
\newcommand{\loglt}{\log_{\sG}}
\newcommand{\fc}{\mathfrak{t}}
\newcommand{\ku}{K_u}
\newcommand{\kup}{\ku'}
\newcommand{\kt}{\tilde{K}}
\newcommand{\sGpf}{\sG(\O_K)^{pf}}
\newcommand{\hgm}{\widehat{\mathbb{G}}_m}
\newcommand{\bE}{\overline{E}}

\newcommand{\bPi}{\overline{\Pi}}
\newcommand{\bPit}{\bPi^{\rm{\scriptscriptstyle temp}}}
\newcommand{\Pit}{\Pi^{\rm{\scriptscriptstyle temp}}}
\renewcommand{\pit}[1]{\Pi^{\scriptscriptstyle temp}_{#1}}
\newcommand{\pitk}[2]{\Pi^{\scriptscriptstyle temp}_{#1;#2}}
\newcommand{\pio}[1]{\pi_1({#1})}

\newcommand{\xan}{X^{an}}
\newcommand{\yan}{Y^{an}}

\newcommand{\sM}{\mathscr{M}}

\togglefalse{draft}
\newcommand{\FF}{\cite{fargues-fontaine}}
\iftoggle{draft}{\pagewiselinenumbers}{\relax}

\begin{abstract}
Let $p$ be a prime number. Let $X/E$ be a geometrically connected, smooth, quasi-projective variety over a finite extension $E/\Q_p$. In this paper I demonstrate the existence of isomorphs of the tempered (and hence also \'etale) fundamental group of $X/E$ which are labeled by distinct arithmetic  holomorphic structures, just as  isomorphs of the fundamental group of a Riemann surface $\Sigma$ may be labeled by Riemann surfaces (i.e. complex holomorphic structures) $\Sigma'$ in the Teichmuller space of $\Sigma$. This is the starting point of the theory elaborated in \cite{joshi-teich,joshi-teich-estimates,joshi-teich-summary-comments, joshi-teich-rosetta} for which this paper is intended as an brief sketch and announcement. Arithmetic holomorphic structures introduced here also provide distinct arithmetic holomorphic structures used by Mochizuki in \iut. Since the question of whether or not there exists of distinct arith. hol. structures in \iut\ was raised in \cite{scholze-stix}, I  include a discussion of \cite{scholze-stix}. See the introduction for additional details.
\end{abstract}

\section{Introduction}
\subsection{} 
Classical Teichmuller Theory of Riemann surfaces provides a collection of labeled \emph{isomorphs} (i.e. isomorphic copies) of the topological (and hence also the \'etale) fundamental group of a Riemann surface; in this labeling  each isomorph of the fundamental group is  labeled  by the quasi-conformal structure i.e. the point of the Teichmuller space giving rise to the isomorph (see  Section~\ref{se:untilts-of-Pi-teichmuller}).  One may even view the classical Teichmuller space as providing an analytic family of isomorphs of the fundamental group of a Riemann surface. Here in Section~\ref{se:untilts-of-Pi} and especially Theorem~\ref{thm:arith-hol-strs},  I show that one can explicitly construct topologically/geometrically distinguishable data which provide  isomorphs of the tempered fundamental group of a geometrically connected, smooth, quasi-projective variety over  $p$-adic fields and mirroring the classical case,   I show  that there exist many analytic (and adic) spaces parameterizing  isomorphs of the tempered fundamental group (i.e. the labeling provided here is by natural, continuous and even analytic parameters). 

By analogy with the classical case of Riemann surfaces, the  existence of such labeled isomorphs in the $p$-adic context (as constructed in this paper) can be understood as indicative of the existence of a $p$-adic Teichmuller Theory--this is the theory of \cite{joshi-teich}, \cite{joshi-teich-estimates}. This natural labeling is achieved (in Section~\ref{se:untilts-of-Pi}) by means of Theorem~\ref{thm:main}, Theorem~\ref{thm:main2}.
The isomorphs provided by Theorem~\ref{thm:main} and Theorem~\ref{thm:main2} arise from untilts of any algebraically closed perfectoid field of characteristic $p>0$ and hence I call these isomorphs \emph{untilts of fundamental groups}, or more precisely \emph{untilts of tempered fundamental groups} (one can use this to define untilts of \'etale fundamental groups as well). Untilts of fundamental groups constructed here all arise from \textit{Arithmetic Holomorphic Structures} (Definition~\ref{def:arith-hol-strs}) and every arithmetic holomorphic structure in the sense of Definition~\ref{def:arith-hol-strs} also provides an arithmetic holomorphic struture in the sense of \iut\ i.e. in the sense of \cite[Example 1.8]{mochizuki-iut2} (see Theorem~\ref{thm:arith-hol-strs}) (see Subsection~\ref{ss:intro-mochizuki-context} for more on this).

Let me also remark that the key observation which was the starting point for this note (and the theory of arithmetic Teichmuller spaces which I have elaborated in \cite{joshi-teich}, \cite{joshi-teich-estimates}, \cite{joshi-teich-summary-comments}, \cite{joshi-teich-rosetta}) was Theorem~\ref{thm:main2}, which also shows that the absolute Grothendieck conjecture fails for the class of Berkovich spaces (over algebraically closed perfectoid fields), arising as analytifications of geometrically connected, smooth, projective variety over  $p$-adic fields (just as the conjecture fails for complex analytic spaces especially for Riemann surfaces). In particular Theorem~\ref{thm:main2} justifies the use of the term arithmetic holomorphic structures in Definition~\ref{def:arith-hol-strs} and in Theorem~\ref{thm:arith-hol-strs}.

\newcommand{\fjxe}{\mathfrak{J}(X/E)}
\newcommand{\fJ}{\mathfrak{J}}
In Section~\ref{se:arith-hol-spaces}, I show how the definition of arithmetic holomorphic structures given in Section~\ref{se:untilts-of-Pi}  provides an anabelian construction of Arithmetic Teichmuller spaces $\fjxe$ where $X/E$ is a geometrically connected, smooth quasi-projective variety over a $p$-adic field. Remark~\ref{re:arithmetic-space-point} demonstrates the highly non-trivial anabelian nature of this space for the special case when $X=\Spec(E)$ (i.e. a point). Theorem~\ref{th:groth-and-connectedness} establishes the relationship between the absolute Grothendieck conjecture over $p$-adic fields and the connectedness of arithmetic Teichmuller spaces. Notably, as the absolute Grothendieck conjecture is known (by Mochizuki's work)  for the case of a smooth, geometrically connected hyperbolic curve of strict Belyi type over a $p$-adic field one deduces that the arithmetic Teichmuller space  is connected in this case.

Significance of these results, and their relationship to \iut\ and \cite{scholze-stix} is explained in the paragraphs  \ref{ss:intro-mochizuki-context}, \ref{ss:Qpbar}, \ref{ss:replacements-iut}, \ref{ss:corollary3-12} %
 (\textcolor{red}{readers not interested in these two topics may skip these subsections completely} and proceed to Section~\ref{se:main}, Section~\ref{se:untilts-of-Pi} for the main theorems of this note). 

The proofs presented here are adequate for expert readers (detailed proofs appear in \cite{joshi-teich,joshi-teich-estimates,joshi-teich-rosetta}) and this note also serves as a \textit{minimal} announcement (i.e. containing briefest sketch) of the arithmetic Teichmuller Theory originating from these two theorems detailed in \cite{joshi-teich,joshi-teich-estimates,joshi-teich-summary-comments,joshi-teich-rosetta}, and its application to the claims of \iut.

\subsection{}\label{ss:intro-mochizuki-context} To underscore the importance of the results of this paper and \cite{joshi-teich,joshi-teich-estimates,joshi-teich-rosetta,joshi-teich-summary-comments}, let me point out the following. In \iut, Mochizuki asserts  the existence of distinct group theoretic data i.e. data of a surjection \be\label{eq:mochizuki-arith-strct} \Pi=\pi_1^{temp}(X/E)\twoheadrightarrow G_E=G\ee of topological groups (such a surjection is called an arithmetic holomorphic structure in \iut--especially \cite[Example 1.8]{mochizuki-iut2} and this is the sense in which the term is used in the principal results of \cite{mochizuki-iut3}). %

To be clear, it is asserted in \cite[\S 1, Page 229, first paragraph]{mochizuki-iut2}, that the theory developed in \iut\ is about comparing one such arithmetic holomorphic structure with respect to another and hence,   one must  have (a) many such structures available, (b) one must be given a way of distinguishing between two such structures, and (c) establish how the key operations of the theory change such structures. To paraphrase \cite{scholze-stix}, these three points are not transparently established in \iut\ (I am in agreement with \cite{scholze-stix} on this point). 

On the other hand, \cite[Remark 9]{scholze-stix} also asserts (because on the usage   of \cite[Theorem 1.9]{mochizuki-topics3}  in \iut\ and  the validity of the absolute Grothendieck conjecture in the context of strict Belyi Type hyperbolic curves used in  \iut) that any such surjection \eqref{eq:mochizuki-arith-strct} of topological groups must arise from geometric data which is unique up to an isomorphism i.e. \cite[Remark 9]{scholze-stix} asserts the rigidity of  geometric data providing isomorphs  of the fundamental group. [\textit{This rigidity claim  is false because of Theorem~\ref{thm:arith-hol-strs}.}] 

[At this juncture I should say that  upon reading the 2020 version of this note \cite{joshi-untilts}, Scholze said in his email (October 2020) to me that I should emphasize here that \textit{\cite{scholze-stix} especially \cite[Remark 9]{scholze-stix}  was written in the context of the claims in \iut\ based on \cite[Theorem 1.9]{mochizuki-topics3} and its variants and their usage in \iut.} Neither this paper nor \cite{joshi-teich,joshi-teich-estimates,joshi-teich-summary-comments, joshi-teich-rosetta} make any use of Anabelian Reconstruction Theory.]

Theorem~\ref{thm:main}, Theorem~\ref{thm:main2} provide a natural  way of defining arithmetic holomorphic structures (Definition~\ref{def:arith-hol-strs}, Lemma~\ref{le:non-isom-arith-hol-strs} and  Theorem~\ref{thm:arith-hol-strs}) and under the hypothesis of Theorem~\ref{thm:main2} these are literally holomorphic structures in the sense of Berkovich's Theory of analytic functions and analytic spaces as arising from deformations of the pair of analytic spaces $(\xan_E,\xan_K)$ arising from deformations of the perfectoid field $K$ i.e. changes to the (Berkovich) holomorphic structures due to topological changes in arithmetic.  \textit{So the term `arithmetic holomorphic structure' really acquires geometric, function-theoretic meaning in the conventional usage  of the term complex holomorphic structures.} Theorem~\ref{thm:arith-hol-strs} also shows that my definition of arithmetic holomorphic structures (Definiton~\ref{def:arith-hol-strs}) also provides (analytic families of) distinct arithmetic holomorphic structures  in Mochizuki's group theoretic sense i.e. a surjection \eqref{eq:mochizuki-arith-strct}. This is the approach of Section~\ref{se:untilts-of-Pi} based on Theorem~\ref{thm:main}, Theorem~\ref{thm:main2}, (Lemma~\ref{le:non-isom-arith-hol-strs} and Theorem~\ref{thm:arith-hol-strs}) and is elaborated in \cite{joshi-teich}.

\textit{Thus my approach (presented in this paper) provides a robust demonstration (Theorem~\ref{thm:arith-hol-strs}) of one of the important unproved (\cite{scholze-stix}) assertions upon which \iut\ ultimately rests and contradicts the rigidity claim of \cite[Remark 9]{scholze-stix} on which \cite{scholze-stix} objections to the strategy of \iut\ ultimately rests. There are other objections raised in \cite{scholze-stix} which also can be dealt with using the results of \cite{joshi-teich}, \cite{joshi-teich-estimates}.}

Notably in \cite{joshi-teich-estimates}, I provide an intrinsic construction of Mochizuki's $\Theta_{gau}$-Link which is central to \cite[Corollary 3.12]{mochizuki-iut3} and is ultimately responsible for a certain $j^2$-factor, appearing in the valuation calculations of \cite[Corollary 3.12]{mochizuki-iut3}  (and hence in \cite{mochizuki-iut4}), whose appearance has been held in doubt by \cite[Section 2.2]{scholze-stix} as Mochizuki does not provide a convincing demonstration of its existence (I am in agreement with \cite{scholze-stix} that the key geometric mechanism which underlies this claim is not convincingly established in \cite{mochizuki-iut3}). 

After  my construction of $\Theta_{gau}$-Links in \cite{joshi-teich-estimates} there is little doubt that this $j^2$ factor does indeed arise for non-trivial reasons and \cite{joshi-teich-estimates} shows that it arises from delicate aspects of Fargues-Fontaine curves (Mochizuki has asserted \cite[Fig. I2, Page 411]{mochizuki-iut3}, presciently, but without proof that $\Theta_{gau}$-Link has to do with Witt vectors--the exact version of this claim is proved in \cite{joshi-teich-estimates}).  

As mentioned earlier, an important consequence of my results  is that there exist many continuous and analytic  families (i.e. families with parameter spaces which are $p$-adic schemes or adic spaces) of isomorphs of the tempered fundamental group of any geometrically connected, smooth, quasi-projective varieties over $p$-adic fields (and hence a similar assertion holds for \'etale fundamental groups). This is an important point not demonstrated in \iut, but it is really needed to serve as a (geometric) firmament for Mochizuki's (group theoretic) view of ``Galois groups as arithmetic tangent bundles'' \cite[Page 24]{mochizuki-iut1} which forms the back drop for \iut. \textit{Importantly  in each of these families of isomorphs constructed here (Theorem~\ref{thm:main}, Theorem~\ref{thm:main2} and Theorem~\ref{thm:arith-hol-strs}), the fixed prime number $p$ behaves as a continuous parameter (by \cite{fargues-fontaine}). This behavior of $p$ has been sought for in many Diophantine applications, including the $abc$-conjecture.}

Mochizuki, in his comments (see \cite{mochizuki-ss-rejoinder}, \cite{mochizuki-essential-logic}) on \cite{scholze-stix} has asserted (see \cite[\S 3]{mochizuki-essential-logic}) that  there are no redundant isomorphs in \iut, while \cite{scholze-stix} have argued, based on assertions in \iut\ and \cite{mochizuki-topics3}, that it is enough to work with one isomorph of the tempered fundamental group given that all of them provide isomorphic geometric data (i.e. all but one isomorphs are redundant).  On the other hand no proof of Mochizuki's non-redundancy claim has appeared to date.  \textit{This paper  provides the first proof of Mochizuki's non-redundancy claim by establishing that the isomorphs are of distinct arithmetic-geometric provenance (and even continuous families of isomorphs exist) and therefore are non-redundant (in both \iut\ and \cite{joshi-teich}--note that the latter works directly with geometric data providing the isomorphs).}

\subsection{}\label{ss:Qpbar} Let me explain how algebraically closed, perfectoid fields are required for \iut\ (but are never considered in it). Algebraically closed, complete (rank one) valued fields equipped with isometric embeddings of $\Q_p$  necessarily enter \iut\ through the critical role which arbitrary geometric base-points play in \iut\ (see \cite[\S I3, Page 26]{mochizuki-iut1} for Mochizuki's discussion of the  role of arbitrary geometric base-points in his theory; my discussion of geometric base-points for tempered fundamental groups is in Section~\ref{se:untilts-of-Pi}). \textit{Such fields are perfectoid by Lemma~\ref{lem:perfectoid} and the field $\C_p$ is just one example.}  Let me remark as an aside that Mochizuki builds his theory (\iut) using the field $\bQ_p$ (which is not complete). On the other hand note that the algebraic closure of $\Q_p$ in an algebraically closed perfectoid field $K$ may not be dense in $K$ even if one fixes the tilt to be $K^\flat\isom\cpt$--this is due to \cite{kedlaya18}.   A detailed examination of \iut\ and \topics, and my conversations with readers and experts familiar loc. cit. suggests to me that these papers tacitly work under the assumption that the (only) relevant field for \iut\ is $\C_p$ (\textit{this is contrary to the assertion \cite[\S I3, Page 26]{mochizuki-iut1} regarding geometric base-points in \iut}). \textit{To be perfectly clear,} Mochizuki assiduously avoids  the usage of such large fields (for example $\C_p$) altogether, because of Anabelian Reconstruction Theory does not reconstruct such fields. I came to recognize in the course of writing \cite{joshi-teich,joshi-teich-estimates} that this discrepancy leads  to many of the issues raised by \cite{scholze-stix} regarding \iut. Once the presence of algebraically closed, perfectoid fields in the theory  is understood, it is possible to establish a natural equivalence between local,  i.e. $p$-adic for each prime, initial data of the theories of \cite{joshi-teich} and \iut\ (this is detailed in \cite{joshi-teich-rosetta}). While \cite{scholze-stix} discovered  inadequacies of \iut, but \cite{scholze-stix} did not precisely recognize the source of these issues. I was the first to recognize the source (through this paper) and my work \cite{joshi-teich,joshi-teich-estimates,joshi-teich-summary-comments,joshi-teich-rosetta} provides a fundamentally original resolution of many of the issues noted by \cite{scholze-stix}. 

I have described the  theory which emerges from this approach in \cite{joshi-teich,joshi-teich-estimates,joshi-teich-summary-comments,joshi-teich-rosetta}. By its very construction, my theory does not suffer from any of the drawbacks of \iut\ highlighted by \cite{scholze-stix} and also establishes many of its claims from an independent perspective, moreover my results can also be applied to \iut, and so my conclusion regarding \iut, in contrast with that of \cite{scholze-stix} (and \cite{scholze-review}), is that the remarkable strategy of \iut\ to prove Diophantine inequalities presently remains a viable one.

\subsection{}\label{ss:replacements-iut} Let me remark that  I arrived at my theory (\cite{joshi-teich,joshi-teich-estimates,joshi-teich-summary-comments}) independently of \iut, and therefore I felt no compelling reason to discuss \cite{scholze-stix,scholze-review} in \cite{joshi-teich,joshi-teich-estimates,joshi-teich-summary-comments}. On the other hand  \cite{joshi-teich-rosetta} includes a `Rosetta Stone' for translating between objects of \cite{joshi-teich,joshi-teich-estimates} and those of \iut, and also a full discussion  of  \cite{scholze-stix,scholze-review} and how its objections to \iut\ can be circumvented by applying results of \cite{joshi-teich}, \cite{joshi-teich-estimates} to \iut. \textit{Here are some examples of the required changes provided by my theory.}

\para The notion of arithmetic holomorphic structures in \iut\ is replaced with the conceptually cleaner definition of arithmetic holomorphic structures provided in Definition~\ref{def:arith-hol-strs}.  Theorem~\ref{thm:arith-hol-strs} then provides distinctly labeled isomorphs (labeled by arithmetic holomorphic structures) of the tempered (and hence \'etale) fundamental groups of hyperbolic curves required for \iut. Notably now isomorphs are labeled by distinct arithmetic holomorphic structures in the sense of Definition~\ref{def:arith-hol-strs} (and this agrees with our conventional notion of holomorphic structures by Theorem~\ref{thm:main2}). By Theorem~\ref{thm:arith-hol-strs}, this also agrees with arithmetic holomorphic structures in the sense of \iut.
\para The fundamental `prime-strip' $G_E\act \bQ_p^{\times\mu}$  of \cite{mochizuki-iut3} (which plays a central role and essentially provides Mochizuki's Indeterminacy Ind2 in \cite[Theorem 3.11, Corollary 3.12]{mochizuki-iut3}), here $E$ is a $p$-adic field and $G_E$ is its absolute Galois group computed using some fixed algebraic closure of $E$, is replaced in my theory \cite{joshi-teich} (especially \cite[Theorem 8.29.1]{joshi-teich}) by $G_E\act \hgm(\O_{\cpt})$ (where $\hgm/\Z_p$ is the multiplicative formal group) i.e.  the fundamental monoid of the theory is not $\bQ_p^{\times\mu}$ (as \iut\ essentially asserts as it is formulated using this monoid) but $\hgm(\O_{\cpt})=(1+\mathfrak{m}_{\cpt})$ (\cite{joshi-teich,joshi-teich-estimates,joshi-teich-rosetta}). [By \cite{fargues-fontaine}, one has a natural identification (of $\Q_p$-Banach spaces) $\hgm(\O_{\cpt})=B^{\varphi=p}$ (compatible with the action of $G_E$) where $B=B_{\cpt,\Q_p}$ is the canonical ring of $p$-adic Hodge Theory. This makes geometric and $p$-adic Hodge theoretic role and meaning of considering  $G_E\act \hgm(\O_{\cpt})$ completely transparent. Each untilt $(K,K^\flat\isom \cpt)$ of $\cpt$ provides an isomorph $\invlim_{x\mapsto x^p}\hgm(\O_K)\isom \hgm(\O_{\cpt})$ of the fundamental monoid and by transport of structure, an isomorph fundamental prime-strip of my theory, namely $G_E\act\left(\invlim_{x\mapsto x^p}\hgm(\O_K)\right)$ etc. making it completely transparent that there are many distinctly labeled isomorphs of the fundamental prime-strip $G_E\act \hgm(\O_{\cpt})$ arising from distinct arithmetic-geometric data (as required for \iut). [At this juncture I should say that the topological group action datum $G_E\act\hgm(\O_{\cpt})$ does not uniquely identify $\mathscr{X}_{\cpt,E}$ \cite{joshi-gconj} and this is one of the origins of Mochizuki's first indeterminacy Ind1.]   In particular this means the Mochizuki's Ism-group ${\rm Aut}(G_E\act \bQ_p^{\times\mu})$ which appears in \cite[Theorem 3.11 and Corollary 3.12]{mochizuki-iut3} is replaced in my theory by ${\rm Aut}(G_E\act \hgm(\O_{\cpt}))$ etc. this makes it possible to provide a clear demonstration of how arithmetic holomoprhic structures change in \cite{joshi-teich,joshi-teich-estimates,joshi-teich-summary-comments,joshi-teich-rosetta} (the existence of such changes is claimed in  \iut). [In this context, another important discovery of \cite{joshi-teich} is that the precise version of Mochizuki's indeterminacy Ind2 which is described above via ${\rm Aut}(\hgm(\O_{\cpt}))$ has a classical analog: namely it is the action of the \textit{Virasoro algebra} (or the \textit{Virasoro group} i.e the group (scheme) denoted by ${\rm Aut}(\mathcal{K})$ in \cite[Proof of Theorem~17.3.2]{frenkel01-book}) on  Teichmuller spaces and moduli of Riemann surfaces see \cite{kontsevich87}, \cite{beilinson88} (or \cite{frenkel01-book} for a modern treatment)!] I do not use the term `prime-strips' or any other terminology of \iut. 

\para Mochizuki's $\Theta_{gau}$-Link of \cite{mochizuki-iut3} is replaced by the construction of  \cite[Section 7]{joshi-teich-estimates} and I provide an intrinsic and explicit construction of the set of such links (this set is called Mochizuki's Ansatz in \cite{joshi-teich-estimates}) which also makes it clear that there are in fact many distinct $\Theta_{gau}$-links etc. Moreover this construction demonstrates the existence of a non-trivial $j^2$ factor alluded to earlier.

\para Other tools used in \iut\ such as Frobenioids \cite{mochizuki-frobenioid1} are not required in my work but can also be dealt with greater precision from my point of view. Here is an example illustrating my claims (taken from \cite{joshi-teich-rosetta}). Let $(E\into K,K^\flat\isom \cpt)$ be an $E$-untilt of $\cpt$ (with $E$ a $p$-adic field), then the Frobenioid of the non-archimedean field $K$ (\cite[Examples 6.1, 6.2, 6.3]{mochizuki-frobenioid1}, \cite[Example 11.4]{fucheng}) is perfect and provides a perfection of the Frobenioid of $E$. If $K$ has value-group $\R$ then the Frobenioid of $K$ even provides a realification of the Frobenioid of $E$ in the sense of \cite{mochizuki-frobenioid1} and \iut. The perfection of a Frobenioid of a $p$-adic field is unique (up to an isomorphism) but as the field $K$ varies in my approach  one sees  at once that Mochizuki's Frobenioidal data (for the perfection of the Frobenioid of $E$) itself varies arithmetically and non-trivially in my theory! Moreover all the value group data provided by $K$ can be read in one common value group given by the value group of the tilt $\cpt$ via the given isomorphism $K^\flat\isom \cpt$ (as was noted in \cite{joshi-teich-estimates}, this is important for \cite[Corollary 3.12]{mochizuki-iut3}). Advantages my approach offers should not be lost on the readers--Mochizuki works with the theory of Frobenioids  all the while requiring  such non-trivial arithmetic (and hence geometric via Theorems~\ref{thm:main}, \ref{thm:main2}) variation in the Frobenioidal (i.e. monoidal) data.

\subsection{}\label{ss:corollary3-12} Finally let me say a few words about  \cite[Corollary 3.12]{mochizuki-iut3} as it has become a focus point in the debate on \iut\ after the appearance of \cite{scholze-stix}. [My approach to this is detailed in \cite{joshi-teich-estimates} and further elaborated in \cite{joshi-teich-rosetta}.]  
The key difficulty in \cite[Corollary 3.12]{mochizuki-iut3} lies in providing the correct construction of the $\Theta$-values set (subject to some rules).

\para First of all let me remark that this corollary has conceptual analogs in the theory of Riemann surfaces  (this is detailed in \cite{joshi-teich-summary-comments,joshi-teich-rosetta}). In the theory of Riemann surfaces, the study of sums (resp. integrals) of positive, real valued functions over (1) Teichmuller spaces, (2)  orbits of the Teichmuller modular groups, or (3) some (extended) moduli spaces is quite common--see \cite{mirzakhani07} or the survey \cite{wright19} for results of this type--notably MacShane's Identity  for $\sM_{1,1}$ (see \cite{wright19}) is a contextually  pertinent example. Mochizuki's \cite[Corollary 3.12]{mochizuki-iut3} must be understood as being of this cadre modulo the construction of arithmetic Teichmuller spaces. Especially, in the $p$-adic context of \cite[Corollary 3.12]{mochizuki-iut3}, the products  of $p$-adic absolute values of the $\Theta$-values (from different arithmetic holomorphic structures) is an example of a real-valued  function on arithmetic Teichmuller spaces--see \cite{joshi-teich-estimates} (in \cite{mochizuki-iut3} one takes the usual logarithm of the $p$-adic absolute values, so the products (\cite{joshi-teich-estimates}) becomes sums). So the existence of such a  corollary  must be viewed as consequent to the existence of (arithmetic) Teichmuller spaces at all primes $p$.

\para The issue highlighted by \cite{scholze-stix} regarding \cite[Corollary 3.12]{mochizuki-iut3} is that  \iut\ is unclear on three points (1) how to demonstrate the existence of distinct arithmetic holomorphic structures? and (2) how exactly does the $\Theta_{gau}$-Link  relate distinct arithmetic holomorphic structures and more importantly how does one demonstrate that it contributes a non-trivial valuation scaling factor of $j^2$? (3) What exactly changes in the column of  $\log$-links used in \cite[Corollary 3.12]{mochizuki-iut3}. Because these three points are not transparently established in \iut\ (according to \cite{scholze-stix}), the existence of the $\Theta$-values locus (central to \cite[Corollary 3.12]{mochizuki-iut3}) itself becomes problematic \cite{scholze-stix}.

\para My work provides a resolution of these issues (my construction of the $\Theta$-values locus is detailed in \cite{joshi-teich-estimates} and elaborated in \cite{joshi-teich-rosetta}).  The  resolution of the first two has been discussed earlier in Subsection~\ref{ss:replacements-iut}. %
The third point is somewhat subtle and \iut\ is rather unclear (to me) on this. This issue is resolved (in my approach, \cite{joshi-teich,joshi-teich-rosetta}) by recognizing that if one works with $\syqp$ (i.e. the adic Fargues-Fontaine curve as opposed to the ``complete curve'' $\sxqp$) then  the ${\bf 1}$-column of $\log$-links in \cite[Corollary 3.12]{mochizuki-iut3} corresponds to the fiber over a point, for example the canonical point $x_{can}\in\sxqp$, of the canonical morphism $\syqp\to \syqp/\varphi^\Z=\sxqp$! The residue fields of all closed classical points in the fiber over $x_{can}\in\sxqp$ are (naturally) isometric \cite[Th\'eore\`me 6.5.2(5)]{fargues-fontaine} and hence arithmetic holomorphic structure in the sense of \iut\ does not change in the fibers of this morphism. This should be compared with Mochizuki's assertion that arithmetic holomorphic structure (in his sense) does not change  in the ${\bf1}$-column of $\log$-links  (see \cite[Proof of Corollary 3.12, Paragraph before Remark 3.12.3]{mochizuki-iut3}).  But arithmetic holomorphic structure in the sense of Definition~\ref{def:arith-hol-strs} does change (in a precisely understood way) and so this makes the relationship completely explicit. 

\para Notably an important consequence of \cite{joshi-teich}, \cite{joshi-teich-estimates} is that one may construct $\Theta$-values loci with values in many different rings of $p$-adic Hodge Theory (this is detailed in \cite{joshi-teich-rosetta}, but \cite{joshi-teich-estimates} should be adequate for experts on $p$-adic Hodge Theory) such as $B$, $B_{cris}$, $B[1/t]$, $B_I \text{ for } I\subset [0,1]\subset \R$, $B_I[1/t]$ and even $B_{cris}, B_{dR,\Q_p}$ (this is not available in \iut). 

\para In \cite{joshi-teich-rosetta}, I also detail a Mochizuki style the construction of the $\Theta$-values locus in Galois cohomology  $H^1(G_E,\Q_p(1))$ more precisely in the Fontaine subspace  $H^1_f(G_E,\Q_p(1))\subset H^1(G_E,\Q_p(1))$. As was observed in \cite[Remark 17.5]{joshi-anabelomorphy}, the Fontaine subspace $H^1_f(G_E,\Q_p(1))$ is Mochizuki's log-shell  tensored with $\Q_p$ (see \cite[\S, Page 407]{mochizuki-iut3} for a discussion of $\log$-shells). From my point of view a better codomain for this Galois cohomology construction is the space $H^1(G_E,B^{\varphi=p})$ which reflects the role of the Fargues-Fontaine curves in the theory.  There is one draw back for the Galois cohomology constructions in \iut: $H^1_f(G_E,\Q_p(1))$ is not a ring but only $\Q_p$-vector space and so it does not make sense to multiply $\Theta$-values  (arising from different arithmetic holomorphic structures and even distinct primes lying over a fixed $p$) in this space. Mochizuki's way of compensating for the lack of multiplication in $H^1_f(G_E,\Q_p(1))$ is the introduction of tensor products of these spaces in \cite[Section 3]{mochizuki-iut3}.  On the other hand $B$ is a ring and multiplication makes sense intrinsically, and I expect that this will simplify the eventual calculations which need to be made.

\interfootnotelinepenalty=10000

\para The results of the present series of papers suggest that it is quite reasonable to expect that by combining Mochizuki's rubric  with the methods of \cite{joshi-teich,joshi-teich-estimates,joshi-teich-summary-comments,joshi-teich-rosetta}, the principal Diophantine assertion of \iut\ can be transparently and independently established.

\subsection*{Acknowledgments} %
\textit{My intellectual debt to Shinichi Mochizuki should be clear.} None of \cite{joshi-formal-groups,joshi-anabelomorphy,joshi-gconj,joshi-teich-estimates,joshi-teich-summary-comments,joshi-teich-rosetta} (and the present paper) could have existed without the backdrop provided by his work on $p$-adic Teichmuller Theory \cite{mochizuki-ordinary}, \cite{mochizuki-teichmuller-book}, anabelian geometry and \iut. Influence of the works of \cite{fargues-fontaine}, \cite{kedlaya18}, \cite{scholze12-perfectoid-ihes} should also be self-evident.

This note is a substantially enhanced and revised version of \cite{joshi-untilts} which is now subsumed in \cite{joshi-teich,joshi-teich-estimates}. While much of the material first presented in \cite{joshi-untilts} has been replaced or enhanced here, that note certainly forms the skeleton of this one  and the original title has been retained (but a subtitle has been added).  That version was put into a limited circulation some time in the Spring of 2020 and was posted on the arxiv subsequently. So this paragraph of the Acknowledgments is borrowed from \cite{joshi-untilts}.  My thanks are due to Peter Scholze,   Emmanuel Lepage, and Jacob Stix, for  promptly (in 2020) providing comments, suggestions or corrections to \cite{joshi-untilts}. I thank Yuichiro Hoshi for answering (in 2020) some of my elementary questions  about tempered fundamental groups in the early days of my investigations. Thanks are due to Kiran Kedlaya for comments on \cite{joshi-untilts} and early versions of \cite{joshi-teich}, \cite{joshi-teich-estimates} and his encouragement. 

It is a great pleasure to thank one  mathematician who wishes to remain anonymous, but who insisted that I provide a full detailed proof of Theorem~\ref{thm:main2} (this detailed proof appears in \cite{joshi-teich}) and made a number of suggestions regarding this proof. I have benefited greatly from his general encouragement and sagely counsel.

\newcommand{\ebh}{\widehat{\bE}}
\newcommand{\ebhx}[1][x]{\widehat{\bE^{#1}}}
\newcommand{\bdr}{B_{dR}}
\newcommand{\bdre}{{\bdr}_{,E}}
\newcommand{\bdrep}{B^+_{dR,E}}
\newcommand{\kbh}{\widehat{\bK}}

\numberwithin{equation}{section}
\section{Untilts of fundamental groups of Riemann surfaces}\label{se:untilts-of-Pi-teichmuller}
Before discussing the $p$-adic theory outlined in  Section~\ref{se:main}, Section \ref{se:untilts-of-Pi}
let me illustrate my approach using the theory of Riemman surfaces. 

Let $\Pi_\Sigma=\pi_1^{top}(\Sigma)$ be the topological fundamental group of a connected Riemann surface $\Sigma$, which one assumes to be hyperbolic to avoid trivialities. Then consider  connected Riemann surfaces $\Sigma'$ whose topological type is the same as that of $\Sigma$, so that for any choice of a homeomorphism $f:\Sigma\to \Sigma'$ one has an isomorphism of topological fundamental groups $\pi_1^{top}(\Sigma')\isom \Pi_\Sigma$ and let me write the information provided by this isomorphism as $\Pi_{\Sigma;\Sigma'}$. I will call $\Pi_{\Sigma;\Sigma'}$ an \textit{untilt of $\Pi_\Sigma$ provided by the holomorphic structure $\Sigma'$} and  in this way the complex structure of $\Sigma'$ serves as a geometrically distinguishable label of this isomorph of $\Pi_\Sigma$. \textit{Thus the labeling of isomorphs of  $\Pi_{\Sigma}$  contemplated here  is by holomorphic structures $\Sigma'$ on the topological space $\Sigma$.}  

[Note that I will suppress from my notation $\Pi_{\Sigma,\Sigma'}$ the choice of the homeomorphism $f$ because, in the context of Teichmuller space of $\Sigma$, such a homeomorphism is included in the datum of any point of the Teichmuller space \cite[Chap. V, Section 2.1]{lehto-book}. Note also that in the classical context I have suppressed the  dependence on the choice of base-points, but base-points  play an important role in $p$-adic context as is discussed in Section~\ref{se:untilts-of-Pi}.]

One can also consider $\Sigma'$ up to some  natural equivalence relations.  For example, let $(g,n)$ be the topological type of $\Sigma$, and consider the classical Teichmuller space $T_{g,n}$ of $\Sigma$ and consider $\Sigma'\in T_\Sigma$, similarly one may also consider the moduli of Riemann surfaces of type $(g,n)$ and consider the isomorphism class of $\Sigma'\in \sM_{g,n}$. The classical group actions on Teichmuller spaces provide actions on the set of isomorphs $\Pi_{\Sigma;\Sigma'}$ (by operating on the geometric label  $\Sigma'$). This leads to the following proposition:

\bthm
Let $\Sigma$ be a connected Riemann surface of topological type $(g,n)$ and let $T_{g,n}$ be the Teichmuller space of $\Sigma$ and let $\sM_{g,n}$ be the moduli of Riemann surfaces. Then 
\benumlab
\item $\Sigma'\in T_{g,n}\mapsto \Pi_{\Sigma;\Sigma'}$ (resp. $\Sigma'\in \sM_{g,n}$) provides a natural continuous family of isomorphs of $\Pi_\Sigma$ labeled by points of the classical Teichmuller space of $T_{g,n}$  (resp. the moduli space $\sM_{g,n}$).
\item Notably the fiber of the natural morphism $T_{g,n}\to \sM_{g,n}$  over the isomorphism class of $[\Sigma]\in\sM_{g,n}$ provides a collection of isomorphs $\Pi_{\Sigma;\Sigma'}$ in which $[\Sigma']=[\Sigma]$ (i.e. the isomorphism class of $\Sigma$ remains fixed).
\item There are two fundamental  actions on the set of isomorphs $\left\{\Pi_{\Sigma,\Sigma'}:\Sigma'\in T_{g,n} \right\}$ given by the action of the following groups on the labels $\Sigma'$:
\benum 
\item the action of the Teichmuller modular group acting on $T_{g,n}$ such that the quotient is $\sM_{g,n}$ (see \cite{lehto-book});
\item and the action of the Virasoro algebra (or the Virasoro group) described in \cite{kontsevich87}, \cite{beilinson88}, \cite[Theorem 17.3.2]{frenkel01-book}.
\eenum
\eenum
\ethm

Let me remind the reader what this means for  $(g,n)=(1,1)$: here $T_{1,1}=\mathfrak{H}$ is the upper half-plane and $T_{1,1}=\mathfrak{H}\to \sM_{1,1}=\mathfrak{H}/{\rm SL}_2(\Z)$ is the (usual) quotient of the upper half-plane by the classical modular group ${\rm SL}_2(\Z)$. \textit{Important observation of the present series of papers is that the above theorem has an arithmetic analog--namely Theorem~\ref{thm:arith-hol-strs} and especially Theorem~\ref{th:main3}.}

\brem 
In Section~\ref{se:main}, the key new idea is to work with deformations of Banach ring structure of functions on varieties over $p$-adic fields. There is in fact a parallel for this in the theory of Riemann surfaces. Notably, classical Teichmuller Theory can be understood as changes of (infinite dimensional) Banach ring structure of suitably nice (but not necessarily holomorphic) functions on a Riemann surface (even in the case of compact Riemann surfaces, one has to work with infinite dimensional Banach spaces). This is discussed in detail in \cite[Section 8]{joshi-teich}. So both the $p$-adic case discussed in Section~\ref{se:main} and classical Teichmuller Theory can be understood as deformations of Banach ring structures of functions!
\erem

Now suppose that $L$ is a number field and $X/L$ is a hyperbolic, geometrically connected smooth quasi-projective curve. By  \cite{tamagawa97-gconj}, \cite{mochizuki96-gconj} the anabelomorphism class of $X/L$ (i.e. the isomorphism class of $\pi_1^{et}(X/L)$ coincides with the isomorphism class of the $\Z$-scheme $X$ (this  assertion is the absolute Grothendieck conjecture as proved in \cite{tamagawa97-gconj}, \cite{mochizuki96-gconj}).
In the terminology introduced in \cite{joshi-anabelomorphy}, the $\Z$-scheme $X$ is \textit{amphoric} i.e. an invariant of the anabelomorphism class of $X/L$. [The terms \textit{anabelomorphisms} and \textit{amphoric} were introduced in \cite{joshi-anabelomorphy}.]

Let $L\into \C$ be any embedding, let $\xan_\C=(X\times_{L\into\C}\C)^{an}$ be the complex analytic space associated to the algebraic curve $X_\C$ (i.e. the associated Riemann surface).  

An \emph{untilt of $\pi_1^{\acute{e}t}(X/L)$ at $\infty$} (here $\infty$ is short for ``\emph{at an archimedean prime $L\into \C$}'') is the data consisting of an embedding $L\into \C$, a  Riemann surface $\Sigma'$ and a quasi-conformal mapping $\Sigma'\to \Sigma$, such that $\Sigma'$ has the same  conformal equivalence class as that of $\Sigma$ i.e. in the moduli of Riemann surfaces of topological type of $\Sigma$ one has  $[\Sigma']=[\Sigma]\in \sM_{g,n}$. One may now define isomorphism of two untilts of $\pi_1^{\acute{e}t}(X)$ at $\infty$  in a natural way. The condition $[\Sigma']=[\Sigma]$ implies that $\Sigma'$ is in the fiber over $[\Sigma]$ of the Teichmuller mapping $T_{g,n}\to \sM_{g,n}$ i.e. the orbit of $\Sigma\in T_{g,n}$ under the Teichmuller modular group. 

\newcommand{\cM}{\mathcal{M}}
Thus one has the following:
\bpro\label{pr:untilts-for-riemann-surfaces}
Let  $X/L$ be a geometrically connected, smooth, hyperbolic curve over a number field $L$ with no real embeddings. Then the isomorphism classes of untilts of $\pi_1^{\acute{e}t}(X/L)$ at $\infty$ are in bijection with $$\widetilde{Hom(L,\C)}\times \left\{\Sigma'\in T_{g,n}: [\Sigma']=[\Sigma]\in \sM_{g,n}\right\},$$ where $\widetilde{Hom(L,\C)}$ is the set of equivalence classes of embeddings of $L \into \C$.
\epro

\brem 
Constructions and the  theorems of this paper  in Section~\ref{se:main} and Section~\ref{se:untilts-of-Pi} are of this type, and exist for all primes and I hope this discussion makes my assertion that the existence of labeled isomorphs of tempered fundamental groups provided here and in \cite{joshi-teich} should be viewed as the existence of arithmetic Teichmuller spaces. This theory is detailed in \cite{joshi-teich,joshi-teich-estimates,joshi-teich-summary-comments,joshi-teich-rosetta}.
\erem

\brem
Owing to the topological rigidity of algebraically closed fields complete with respect to an archimedean absolute value, forced by Ostrowski's Theorem (see Remark \ref{re:archimedean-case}), one could say that  untilts of topological fundamental groups at $\infty$ (i.e. at archimedean primes) can  arise only from the existence of geometric variations of the underlying objects. 
\erem

\brem 
Readers familiar with the classical Szpiro inequality (for surfaces fibered over curves) and its several  different proofs, may notice that the above proposition provides a unified way of viewing these proofs  as taking place  over the ``space of untilts.'' More precisely the ``space of untilts'' provides the geometric Kodaira-Spencer classes which underly these proofs. 
\erem

\section{The main theorem}\label{se:main}
\emph{All valuations on base fields considered in this paper will be rank one valuations.} For the theory of tempered fundamental groups see \cite{andre-book,andre03} or \cite{lepage-thesis}. Berkovich spaces will be strictly analytic (and mostly will arise as analytifications of geometrically connected smooth quasi-projective varieties).

In what follows I will work with algebraically closed, perfectoid fields \cite[Definition 3.1]{scholze12-perfectoid-ihes}. Such fields can also be characterized in many different ways (see \cite[Lemma 6.1.9]{scholze-weinstein}).  For readers unfamiliar with perfectoid fields, the following lemma provides a translation of this condition into  more familiar hypothesis.

\blem\label{lem:perfectoid}\ 
\benumlab
\item Any algebraically closed, complete, non-archimedean valued field (rank one valuation) of characteristic $p>0$ is perfectoid.
\item Any algebraically closed, complete non-archimedean valued field (rank one valuation) equipped with an isometric embedding of $\Q_p$ is an  algebraically closed, perfectoid field (of characteristic zero). 
\eenum
\elem 
\bp 
Both the assertions are immediate from \cite[Definition 3.1]{scholze12-perfectoid-ihes}.  Let me provide a proof using another useful characterization of perfectoid fields \cite[Lemma 6.1.9]{scholze-weinstein}. Let $K$ denote a field satisfying hypothesis {\bf(1)} or {\bf(2)}. Let $R=\O_K$ be the valuation ring of $K$ which provides the given valuation of $K$. It suffices to verify that three hypothesis of \cite[Lemma 6.1.9]{scholze-weinstein-book} are satisfied. As $K$ is algebraically closed, the valuation of $K$ is not discrete, secondly the requirement $\abs{p}_K<1$ holds trivially for {\bf(1)} as $p=0$ and for {\bf(2)} as $\Q_p\into K$ is an isometric embedding. So it remains to verify the third hypothesis of \cite[Lemma 6.1.9]{scholze-weinstein-book} that the Frobenius $\phi:R/p \to R/p$ is surjective. If $K$ has characteristic $p>0$, then obviously  $R=R/p$ (as $p=0$) and $\phi$ is obviously surjective on  $R$, as $K$ (and hence  $R$) is perfect. This proves {\bf(1)}. Now suppose $K$ has characteristic zero. Let $\bar{x}\in R/pR$ and suppose $x\in R$ is an arbitrary lift of $\bar{x}$. Then as $K$ is algebraically closed field, choose $z$ so that $z^p=x$ then $z\mod{pR}$ satisfies $\phi(z\mod{pR})=\bar{x}$ and hence this establishes that $\phi$ is surjective on $R/p$ and proves {\bf(2)} and completes the proof of the lemma.
\ep

For an algebraically closed, perfectoid field $K$ as above, one has naturally associated field, denoted $K^\flat$, which is an algebraically closed, perfectoid of characteristic $p>0$, called the \emph{tilt of $K$} (see \cite[Lemma 3.4]{scholze12-perfectoid-ihes}). If $K$ has characteristic $p>0$ then $K^\flat=K$.

\brem Here is an important example of an algebraically closed perfectoid field. Let $\bQ_p$ be an algebraic closure of $\Q_p$ and let $\C_p$ be the completion of $\bQ_p$. By Lemma~\ref{lem:perfectoid}, $\C_p$ is an algebraically closed perfectoid field. The tilt, $\cpt$ of $\C_p$, can be explicitly described as $$\cpt=\widehat{\overline{\F_p((t))}}$$ for some indeterminate $t$ (see \cite[Example 2.1]{morrow-bourbaki}). By \cite[Definition]{scholze12-perfectoid-ihes} or Lemma~\ref{lem:perfectoid}, $\cpt$ is an algebraically closed, perfectoid field of characteristic $p>0$.
\erem

Let $E$ be a $p$-adic field. Let $F$ be an algebraically closed, perfectoid field of characteristic $p>0$. 
Fix an  algebraically closed  field, perfectoid $F$ of characteristic $p>0$. For example readers can assume, that $F= \C_p^\flat$, but it is important to note that in the theory developed in this paper and \cite{joshi-teich}  one allows $F$ to be arbitrary and Definition~\ref{def:arith-hol-strs} also allows $F$ to vary.
By an $E$-\emph{untilt} of $F$, I will mean data $(E\into K,K^\flat\isom F)$ consisting of a  perfectoid field $K$, an isometric embedding $E\into K$ (so $K$ has characteristic zero), and an isometry $K^\flat\isom F$ \cite{scholze12-perfectoid-ihes}, \cite{fargues-fontaine}. If $E=\Q_p$, I will simply refer to $\Q_p$-untilts of $F$ as \textit{untilts} of $F$. Note that as $F$ is algebraically closed, by \cite[Proposition 3.8]{scholze12-perfectoid-ihes}  $K$ is also algebraically closed. Given the  pair $(F,E)$,  by the theory of \cite{fargues-fontaine}, $E$-untilts  of $F$ are parameterized by the  Fargues-Fontaine curve $\syfe$. In particular untilts of $F$ are parameterized by the curve $\syqp$. One also has ``complete'' Fargues-Fontaine curves $\sxfe$ and $\sxqp$ constructed in \cite{fargues-fontaine}. 

Some words of caution about my notation are warranted here: firstly, neither of the curves $\syfe,\syqp,\sxfe,\sxqp$ are of finite type over $\Q_p$ secondly  $\syfe$ and $\syqp$ are adic curves while $\sxfe$, $\sxqp$ can be defined as schemes over $\Q_p$ and hence there is an adic curve associated to both, but I will write the adic and the schematic version with the same notation. This is not conventional in the existing literature but hopefully not cause any confusion in this text. The adic curve $\sxfe$ is the quotient of $\syfe$ by the infinite discrete group acting freely on $\syfe$ and generated by the Frobenius morphism of $\syfe$ and hence one has a natural quotient morphism of adic Fargues-Fontaine curves $$\syfe\to\syfe/\varphi^\Z=\sxfe.$$

Crucial point for this paper (and \cite{joshi-teich}, \cite{joshi-teich-estimates}) is that \emph{there exist untilts of $\cpt$ which provide algebraically closed perfectoid fields of characteristic zero which are not topologically isomorphic.} This is the main result of \cite[Theorem~1.3]{kedlaya18}. Note that all characteristic zero untilts of $\cpt$ have the cardinality of $\C_p$ and are complete and algebraically closed fields and hence are abstractly isomorphic fields but may not be topologically isomorphic after \cite[Theorem 1.3]{kedlaya18}.

Now fix a geometrically connected, smooth quasi-projective variety $X/E$, with $E$ a $p$-adic field. Let $\xan/E$ be the strictly analytic space associated to $X/E$. Let \be\pit{X/E}=\pi_1^{temp}(\xan/E)\ee be the tempered fundamental group of the strictly $E$-analytic space associated to $X/E$ in the sense of \cite{andre03} or \cite{andre-book}.

\emph{Note that my notation $\pit{X/E}$ suppresses the passage to the   analytification $\xan/E$  for simplicity of notation. The theory of (tempered) fundamental groups also requires a choice of geometric base-points which will be suppressed from my notation for the moment, but see the discussion of this in Section~\ref{se:untilts-of-Pi}. The theory of tempered fundamental groups is independent of the choice of geometric base-points by \cite{andre-book}.}

Let $E'/E$ be a finite extension of $E$ with a continuous embedding $E'\into K$ (as $K$ is algebraically closed, valued field containing $E$, such $E'$ exists.   One can consider $X_{E'}=X\times_{E}E'$ (similarly $X_K=X\times_EK$). Then one has, whenever $E'/E$ is Galois, an exact sequence by \cite[Prop. 2.1.8]{andre-book}
$$1\to \pit{X_{E'}/E'}\to \pit{X/E}\to \gal(E'/E)\to 1.$$

Let $\bE\subseteq K$ be the algebraic closure of $E$ contained in $K$.

By  varying $E'$ over all finite extensions of $E\into K$  one obtains (see \cite[Section 5.1]{andre03}) an exact sequence of topological groups:
$$1\to  \ilim_{E'/E}\pit{X_{E'}/E'}\to\pit{X/E} \to \gal(\bE/E)\to 1.$$

\newcommand{\Et}{\widetilde{E}}

\bthm\label{thm:main}
Let $F$ be an algebraically closed perfectoid field of characteristic $p>0$ (for example $F=\cpt$). Let $E$ be a $p$-adic field. Let $K,K_1,K_2$ be algebraically closed perfectoid fields provided by arbitrary untilts of $F$ with continuous embedding $E\into K$ (resp. into $K_1$ and $K_2$). Let $\bE$ (resp. $\bE_1,\bE_2$) be the algebraic closure of $E$ in $K$ (resp. in $K_1,K_2$). Let $X/E$ be a geometrically connected, smooth, quasi-projective variety over $E$. Then one has the following:
\benumlab
\item\label{thm:main-1} a continuous isomorphism $$\pit{X/K}\isom \ilim_{E'/E}\pit{X/E'},$$
where the inverse limit is over all finite extensions $E'$ of $E$ contained in $K$, and a
\item\label{thm:main-2} a short exact sequence of topological groups $$1\to \pit{X/K}\to \pit{X/E}\to G_E\to 1,$$ and
\item\label{thm:main-4} In particular for any two untilts $K_1,K_2$ of $F$, one has a continuous isomorphism $$\pit{X/K_1}\isom \pit{X/K_2}.$$
\eenum
\ethm
\bp 
The assertion \ref{thm:main-1} is the analogue of \cite[Prop. 5.1.1]{andre03} for an arbitrary untilt $K$ of $F$ containing $E$ (as above). Let me remind the reader that my hypothesis on  $K,K_1,K_2$ imply that $K,K_1,K_2$ are algebraically closed and complete with respect to a rank one valuation.  

Let me prove \ref{thm:main-1}, this will also lead to \ref{thm:main-2}. %
Since $K$ is  algebraically closed, it follows that $K$ contains an algebraic closure $\bE$ of $E$.  Let $\Et\subseteq K$ be the closure  (with respect to valuation topology of $K$) of $\bE$. 

It is clear that $\Et\supset \bE$ is complete and algebraically closed field and $\Et$ contains the algebraic closure $\bE\subset K$ of $E$ contained in $K$ with $\bE\subset \Et$ is a dense inclusion. In particular $\Et$ is the completion of $\bE$ with respect to the induced valuation. In other words $\Et$ is a copy of the completion of an algebraic closure of $E$ (usually denoted $\ebh$) equipped with an isometric embedding $\iota:\ebh\into K$ with $\iota(\ebh)=\Et$.
Hence  $K/\Et$ is an isometric extension of algebraically closed, complete valued fields (with rank one valuations). 

Now one can apply  the principle of invariance  of  fundamental groups under passage to extensions of algebraically closed fields. This principle is well-known for \'etale fundamental groups (see \cite[Expos\'e X, Corollaire 1.8]{grothendieck1971a}). For  tempered fundamental groups this principle is proved in \cite[Proposition 2.3.2]{lepage-thesis} (\textit{notably, in \cite{lepage-thesis}, $X/E$ is not required to be proper}). Thus  applying \cite[Proposition 2.3.2]{lepage-thesis} to the extension $K/\Et$ one has an isomorphism of topological groups $$\pit{X_K}\isom \pit{X_{\Et}}.$$

On the other hand by \cite[Proposition 5.1.1]{andre03}, as $\Et$ is the completion of the algebraic closure of $\bE\subset K$ of $E$, one has an isomorphism
\be \pit{X/\Et}\isom \ilim_{E'/E}\pit{X_{E'}/E'}
\ee
and  an exact sequence of topological groups 
\be\label{eq:untilt-sequence} 1\to \pit{X/K} \isom \ilim_{ E'/E}\pit{X_{E'}/E'}\to\pit{X/E} \to \gal(\bE/E)\to 1.\ee
 This proves the assertions \ref{thm:main-1}, \ref{thm:main-2} as claimed.
 
Let me now prove \ref{thm:main-4}. The claimed isomorphism $\pit{X/K_1}\isom \pit{X/K_2}$ follows from the fact that both the groups can be identified with $\ilim_{E'/E}\pit{X/E'}$ where the inverse limit is over all finite extensions of $E'/E$ contained in $K_1$ (resp. $K_2$) and the fact that there is an equivalence between categories of finite extensions of $E$ contained in $K_1$ and the category of finite extensions of $E$ contained in $K_2$, since finite extensions of $E$ are given by adjoining roots of polynomials with coefficients in $E$ and this data is independent of the embedding of $E$ in $K_1$ or $K_2$ and moreover any abstract isomorphism of discretely valued fields is in fact an isometry--i.e given a finite extension of $E$,   $E'\into K_1$ contained in $K_1$, there is an isometry $E'\into K_2$ and vice versa.
\ep

\bthm\label{thm:main2}
Let $X/E$ be a geometrically connected, smooth projective variety. Let $K_1,K_2$ be two untilts of $\cpt$ which contain $E$. Suppose that $K_1,K_2$ are not topologically isomorphic. Then
\benumlab
\item\label{thm:main2-1} one has  an isomorphism of topological groups $$\pit{X/K_1}\isom \pit{X/K_2},$$
\item\label{thm:main2-2} but $\xan/K_1$ and $\xan/K_2$ are not isomorphic $\Q_p$-analytic Berkovich spaces.
\item\label{thm:main2-3} In particular the absolute Grothendieck conjecture fails in the category of Berkovich spaces over perfectoid fields of characteristic zero.
\eenum
\ethm
\bp 
After Theorem~\ref{thm:main},
only \ref{thm:main2-2}  needs to be proved as \ref{thm:main2-2} $\implies$ \ref{thm:main2-3}. The hypothesis of Theorem~\ref{thm:main2} are non-vacuous--by\cite{kedlaya18}, fields $K_1,K_2$ exist. 

Assume that $X/E, K_1,K_2$ are as in my hypothesis and that $X$ is geometrically connected, smooth and projective over $E$. Suppose, if possible, that $\xan/K_1$ and $\xan/K_2$ are isomorphic as strictly analytic Berkovich spaces. Then one has a bounded isomorphism of Banach rings $$K_1\isom H^0(\xan/K_1,\O_{\xan/K_1}) \isom H^0(\xan/K_2,\O_{\xan/K_2})\isom K_2.$$ This isomorphism evidently extends  to a bounded  isomorphism of Banach fields $K_1\isom K_2$ and hence one sees that $K_1$ and $K_2$ are topologically isomorphic. Thus one has arrived at a contradiction.
\ep

\brem 
A substantially detailed proof of \ref{thm:main2-2} appears in \cite{joshi-teich}.
\erem

\brem
As an aside let me remark that the proof of \cite{kedlaya18} (also see \cite[Th\'eor\`eme 2 and \S3 Remarque 2]{matignon84}) provides an uncountable collection of perfectoid fields $K_1,K_2$ with tilts isometric to $\cpt$ and such that  $K_1,K_2$ are not topologically isomorphic.
\erem

\section{Arithmetic holomorphic structures and Untilts of tempered fundamental groups}\label{se:untilts-of-Pi}
The results of Section~\ref{se:main} can be applied to the problem of producing \emph{labeled copies of the tempered fundamental groups in a manner similar to the classical case described in Section~\ref{se:untilts-of-Pi-teichmuller}.}

Let $\sM(A)$ be the Berkovich spectrum \cite{berkovich-book} of any Banach algebra $A$, and let $\xan_E,\xan_K$ be the analytic spaces arising from $X/E$ etc. (as before $X/E$ is a geometrically connected, smooth, quasi-projective variety over a $p$-adic field $E$). It will be useful to recall the definition of geometric base-points in the context of \cite{andre-book}. 

A \textit{geometric base-point} (more precisely a \textit{$K$-geometric base-point}) of the analytic space $\xan_E$   is an algebraically closed, complete (rank one) valued field $K$ and an isometric extension  of valued fields $K/E$ and a morphism (of analytic spaces) $*_K:\sM(K)\to \xan_E$ (see \cite[Chapter III, 1.2.2]{andre-book}). [In the context of  \cite{andre-book} which works in the context of  \cite{berkovich-book}, any geometric base-point provides a valued field whose valuation is of rank one.]  

Note especially that a $K$-geometric base-point  provides an algebraically closed, complete (rank one) valued field $K$ containing $\Q_p$ isometrically and  so by Lemma~\ref{lem:perfectoid} a geometric base-point provides  an algebraically closed, perfectoid  field $K\supset E\supset \Q_p$. The discussion of \cite[\S I3, Page 21]{mochizuki-iut1} makes it clear that \iut\ requires arbitrary geometric base-points (and so the valued field $K\supset E$ must be arbitrary). \textit{This is how algebraically closed, perfectoid fields enter \iut\     (see \cite{joshi-teich-rosetta}).}   

Let me remark importantly that by \cite{kedlaya18} that there exist  algebraically closed, perfectoid fields with an isometrically embedded $\Q_p$ and with tilts isometric to $\cpt$, but such that $\bQ_p$ is not necessarily dense in these fields. So one cannot hope to work (solely) over $\bQ_p$ once one recognizes that arbitrary algebraically closed perfectoid fields are required in \iut\ via its requirement of arbitrary geometric base-points. 

Theorem~\ref{thm:main} and Theorem~\ref{thm:main2} provide us with a way of defining arithmetic holomorphic structures (\cite{joshi-teich}, \cite{joshi-teich-estimates}, especially \cite{joshi-teich-rosetta}) which retains anabelian features as well as our idea of holomorphic functions. [I call it an arithmetic holomorphic structure because Theorem~\ref{thm:main2}(2) shows that  it is a holomorphic structure in the sense of holomorphic functions, secondly this holomorphic structure   is dependent on the arithmetic and topological properties of the coefficient field.] This is the starting point of the Teichmuller Theory presented in \cite{joshi-teich,joshi-teich-estimates} and is elaborated in \cite{joshi-teich-rosetta} (comparison of with \iut\ is in \cite{joshi-teich-summary-comments}).

\begin{defn}\label{def:arith-hol-strs}
Let $X/E$ be a geometrically connected, smooth, quasi-projective variety over a $p$-adic field $E$.	
A (pointed) \textit{arithmetic holomorphic structure} on $X/E$ is the choice of an untilt $(K\supset E,K^\flat\isom F)$ for some algebraically closed perfectoid field $F$ of characteristic $p>0$ and a choice of a $K$-geometric base-point $*_K:\sM(K)\to\xan_E$. \textit{Morphisms of arithmetic holomorphic structures} are be defined as morphisms of the data $((K\supset E,K^\flat\isom F),*_K:\sM(K)\to\xan_E)$ of the arithmetic holomorphic structure.
\end{defn}

\blem\label{le:non-isom-arith-hol-strs}
Non isomorphic arithmetic holomorphic structures exist--more precisely for every perfectoid field $F$ of characteristic $p>0$, there exist arithmetic holomorphic structures $((K_1\supset E,K_1^\flat\isom F),*_{K_1}:\sM(K_1)\to\xan_E)$, $((K_2\supset E,K_2^\flat\isom F),*_{K_2}:\sM(K_2)\to\xan_E)$ which are not isomorphic.
\elem
\bp
This is an immediate consequence of \cite{fargues-fontaine}. By \cite[Corollaire 2.2.22]{fargues-fontaine}, one sees that for any given perfectoid field $F$ of characteristic $p>0$ and any given $p$-adic field $E$, there exist many non-isomorphic untilts and hence in particular there exist many non-isomorphic arithmetic holomorphic structures.
\ep

The choice an arithmetic holomorphic structure also provides the pair of spaces $(\xan_E,\xan_K)$ and the morphism $\xan_K\to\xan_E$ given by base change via $E\into K$. Suppose one is given any $K$-geometric base point $*_K:\sM(K)\to\xan_K$, then
by composition with the base-change morphism $\xan_K\to\xan_E$ provided by the given isometric embedding $E\into K$, one obtains  a $K$-geometric base-point of $*_K:\sM(K)\to \xan_E$ (note the conflation of notation). In particular one can work instead with the datum $((K\supset E,K^\flat\isom F), *_K:\sM(K)\to \xan_K)$ instead of $((K\supset E,K^\flat\isom F), *_K:\sM(K)\to \xan_E)$.

\brem Let me remark (this is not used in what follows, but experts may find it useful) an important observation from the point of view of anabelian considerations  is that arithmetic holomorphic structures on $X/E$ form a category in the obvious way with objects $((K\supset E,K^\flat\isom F),*_K:\sM(K)\to \xan_E)$ where $F$ is arbitrary algebraically closed, perfectoid field of characteristic $p>0$; morphisms between such objects are defined in the obvious way. This category comes equipped with functors to many important categories. For example one has a natural functor to the category of algebraically closed perfectoid fields $$((K\supset E,K^\flat\isom F),*_K:\sM(K)\to \xan_E)\mapsto F.$$ 
As is observed \cite{joshi-teich-rosetta}, the category of (pointed) arithmetic holomorphic structures on $X/E$ is quite closely related to the diamond  $X^\diamond$ associated to the adic space $X^{ad}$  in the sense of the Theory of Diamonds laid out it in \cite[Section 15]{scholze-diamonds}.
\erem

Now I am ready to define untilts of tempered fundamental groups.
\begin{defn}\label{def:untilts-def} Let $X/E$ be a geometrically connected, smooth, quasi-projective variety over a $p$-adic field $E$.
An \emph{untilt of the tempered fundamental group $\Pi=\pit{X/E}$} is the isomorph 
\be\label{eq:untilt-isomorph} \pit{X/E; (K\supset E,K^\flat\isom F)}=\pi_1^{temp}(\xan_E,*_K: \sM(K)\to\xan_E).\ee
of the tempered fundamental group $\pit{X/E}$ computed using the geometric base-point provided by the choice of an arithmetic holomorphic structure $((K\supset E,K^\flat\isom F), *_K:\xan_K\to \xan_E)$ on $X/E$.
\end{defn}

Now suppose $(K_1\supset E,K_1^\flat\isom F), (K_2\supset E,K_2^\flat\isom F)$ are two untilts of $F$ providing arithmetic holomorphic structures on $X/E$, then by \cite{andre-book} one has an isomorphism  of  two topological groups $$\pit{X/E; (K_1\supset E,K_1^\flat\isom F)}\isom \pit{X/E; (K_2\supset E,K_2^\flat\isom F)}.$$ Since untilts may be compared in a natural way (\cite{scholze12-perfectoid-ihes}, \cite{fargues-fontaine}), and hence arithmetic holomorphic structures can be compared. This means that  the geometric and arithmetic data required to define tempered fundamental groups can be compared directly and so arithmetic holomorphic structures (as in Definition~\ref{def:arith-hol-strs}) serve as  natural  geometric labels for isomorphs of tempered fundamental group of $X/E$. 

Theorem~\ref{thm:main2} allows one to assert with impunity that the notion of arithmetic holomorphic structures given in Definiton~\ref{def:arith-hol-strs} is fully compatible with our conventional way of  thinking about distinct classical  holomorphic structures of Riemann surfaces.  As I note in Theorem~\ref{thm:arith-hol-strs}, my definition of arithmetic holomorphic structures also provides arithmetic holomorphic structures in sense this term is used in \iut.

Thus my  approach provides a robust algebro-geometric definition of the phrase arithmetic holomorphic structures and provides a quantitative way of distinguishing such structures. So the sort of problem which has been highlighted by \cite{scholze-stix} about \iut\ does not arise in the theory of the present paper and elaborated in \cite{joshi-teich,joshi-teich-estimates,joshi-teich-summary-comments,joshi-teich-rosetta}.

\brem I will habitually contract the notation $\pit{X/E; (K\supset E,K^\flat\isom F)}$ to $\pit{X/E; K}$ for brevity and simplicity of notation. This may give the impression that the tilting data (i.e. the isometry $K^\flat\isom F$) is unimportant, but let me caution the reader that this is certainly not the case. The tilting data provides  complementary, but necessary, arithmetic information for the valued field $K$ (tilting is ``Frobenius-like'' in the parlance of \iut\ and in my case the term ``Frobenius-like'' is quite literal as  the construction of $K^\flat$ requires $p^{th}$-powers i.e. Frobenius) information in the context of \iut). Notably the datum $K^\flat\isom F$ provided by an untilt allows us to compute arithmetic degrees in one fixed location namely the value group of $F$ even as the untilt moves (see \cite{joshi-teich-estimates} for how this gets used and its relevance for \cite[Corollary 3.12]{mochizuki-iut3}). On the other hand the perfectoid field $K$ provides \'etale like information in the context (and terminology of \iut)--namely it provides the ``\'etale-like'' information given by  the isomorphs $\pit{X/E;K}$ and $G_{E;K}$ of the topological groups $\pit{X/E}$ and $G_E$ respectively.
\erem

By Theorem~\ref{thm:main} an untilt $(K,K^\flat\isom F)$ (for some $F$) and a choice of a geometric base point $\sM(K)\to \xan_E$ provides the short exact sequence of topological groups 
\be\label{eq:galois-seq-untilt} 1\to \pit{X/K}\to \pit{X/E;K}\to G_{E;K}\to 1,\ee	
where the absolute Galois group $G_{E;K}=\gal(\bE_K/E)$ where $\bE_K\subset K$ is the algebraic closure of $E\subset K$ contained in $K$ (the notation $G_{E;K}$  makes the dependence of $G_E$ on the choice of the untilt and the field $K$ it provides   explicit.

\newcommand{\onto}{\twoheadrightarrow}
In particular an arithmetic holomorphic structure in the sense of Definition~\ref{def:arith-hol-strs} also provides arithmetic holomorphic structures in the sense required in \iut\ and one even has the existence of continuous families parameterizing arithmetic holomorphic structures in Mochizuki's sense \cite[Example 1.8, Page 247]{mochizuki-iut2}. By Lemma~\ref{le:non-isom-arith-hol-strs}, non-isomorphic arithmetic holomorphic structures exist. This discussion is summarized in the theorem below.

\bthm\label{thm:arith-hol-strs}
Let $X/E$ be a geometrically connected, smooth, quasi-projective variety over a $p$-adic field $E$. Let $\C_p$ be the completion of a fixed algebraic closure of $\Q_p$. Let $\pit{X/E}$ be the tempered fundamental group of $X/E$ computed using any geometric base-point $\sM(\C_p)\to \xan_E$. Then 
\benumlab
\item non-isomorphic arithmetic holomorphic structures on $X/E$ (Definition~\ref{def:arith-hol-strs}) exist (by Lemma~\ref{le:non-isom-arith-hol-strs}), and 
\item any choice of  a (pointed) arithmetic holomorphic structure   on $X/E$ provides an isomorph  of $\pit{X/E}$ \eqref{eq:untilt-isomorph}, and an exact sequence \eqref{eq:galois-seq-untilt}, and 
\item if $F$ is an algebraically closed, perfectoid field of characteristic $p>0$, then one has continuous families of isomorphs of $\pit{X/E}$ parameterized by Fargues-Fontaine curves $\syfe$ and $\sxfe$.
\item Any choice of an arithmetic holomorphic structure (in the sense of Definition~\ref{def:arith-hol-strs}) provides an arithmetic holomorphic structure   in the sense of \cite[Example 1.8]{mochizuki-iut2}:
$$\pit{X/E;K}\onto G_{E;K}$$
of i.e. an isomorph of the  surjection \eqref{eq:mochizuki-arith-strct} $\pi_1^{temp}(\xan_E)\onto G_E$. 
\item In particular one has continuous families of arithmetic holomorphic structures (in the sense of Definition~\ref{def:arith-hol-strs}) and also in the sense of \iut\ which are parameterized by Fargues-Fontaine curves $\syfe$ and $\sxfe$ for any arbitrary algebraically closed perfectoid field $F$ of characteristic $p>0$.
\eenum
\ethm

Note that the labeling also provides a algebraic substructure of $\Pi$ namely the normal subgroup $\pit{X/K}\subset \pit{X/E}$.

An important discovery of \cite{joshi-teich} is the following (compare with the classical case provided in Proposition~\ref{pr:untilts-for-riemann-surfaces}):
\bthm\label{th:main3} 
Now suppose $F=\cpt$. Then there are two fundamental group actions (acting via their actions on the labels provided by the untilts) on the set $$\left\{ \pit{X/E; (K\supset E,K^\flat\isom \cpt)} : (K\supset E,K^\flat\isom \cpt) \text{ treated as  a } \Q_p\text{-untilt of } \cpt \right\}$$
\benumlab
\item the action of the absolute Galois group $G_E$ on untilts of $\cpt$, (see \cite{fargues-fontaine}) and
\item the action of the group ${\rm Aut}_{\Z_p}(\hgm(\O_{\cpt}))$ (see \cite[Theorem 8.29.1]{joshi-teich}) (which plays the role of the Virasoro group in the classical case described in  Proposition~\ref{pr:untilts-for-riemann-surfaces}).
\eenum
\ethm

\brem \ 
\benumlab
\item As the existence of deformations of complex structures on Riemann surfaces goes hand-in-hand with the existence of non-trivial Kodaira-Spencer classes, similarly the existence of Berkovich analytic structures provided by Theorem~\ref{thm:main2} must be seen as indicative of the existence of Kodaira-Spencer classes of arithmetic interest at every prime $p$. 

\item Since the arithmetic holomorphic structures provided here are closely related to Fargues-Fontaine curves  it is quite reasonable to expect one has non-trivial Kodaira-Spencer classes (as a consequence of Theorem~\ref{thm:arith-hol-strs}, Theorem~\ref{thm:main2}). Examining \cite{mochizuki-HA}, it is quite clear to me that   the Kodaira-Spencer classes  constructed in \cite[Chap IX, \S2]{mochizuki-HA} (for $(g,n)=(1,1)$) must essentially be the same as the ones implicitly provided by Theorem~\ref{thm:main2} (this will be taken up elsewhere). 
\item \textcolor{red}{This remark is not required in this paper or \cite{joshi-teich, joshi-teich-estimates}.} Comparison  \cite[Chap IX, \S2, Page 358]{mochizuki-HA} and especially the assertion $\frac{d}{d\log(q)}\mapsto \left(\frac{d}{dx}\right)^2$  suggests that there is an ``arithmetic heat equation'' (for the classical case see \cite[Chap 1, \S 2]{mumford-theta1}) satisfied by a $p$-adic theta function mirroring the classical case in which the theta-function is annihilated by the differential operator $\frac{d}{d\log(q)}-\left(\frac{d}{dx}\right)^2$ (with  $q=e^{-\pi t}$ in the classical case). Here, in the arithmetic setting, $q$ is regarded as a variable on the arithmetic Teichmuller space constructed \cite{joshi-teich} ($q$ is the Tate parameter i.e. period of a differential of a fixed elliptic curve  over a $p$-adic field). This analogy suggests that one should expect that the classical relationship between periods of differentials on curves and the KP hierarchy to have parallels in the setting of arithmetic Tecichmuller spaces.
\eenum
\erem

\brem
The view of \cite{joshi-teich,joshi-teich-estimates,joshi-teich-summary-comments} is that 
\benumlab
\item Each untilt of $\cpt$ provides a perfectoid field $K$ which is complete and algebraically closed and contains an isometrically embedded $\Q_p$ and hence it provides its own private copy of $\bQ_p\subset K$ i.e. a private copy of $p$-adic arithmetic which may even be regarded as being topologically distinct to the one provided by another such field $K'$ belonging to a different topological isomorphism class from $K$; and 
\item one should view the exact sequence \eqref{eq:galois-seq-untilt} as rising from the pair of analytic spaces $(\xan_E,\xan_K, \xan_K\to\xan_E)$ where the morphism of analytic spaces $\xan_K\to\xan_E$ is induced by base-change using the morphism $E\into K$ of valued fields; 
\item The key discovery of \cite{joshi-teich-estimates} is that Mochizuki's  $\Theta_{gau}$-Link of \cite{mochizuki-iut3} (required for \cite[Corollary 3.12]{mochizuki-iut3}) in fact exists and arises from the geometry of a suitable Fargues-Fontaine curve.
\item (this is the view of \cite{joshi-teich-estimates}) each untilt $(K,K^\flat\isom \cpt)$ of $\cpt$ provides (via the given isometry $K^\flat\isom \cpt$), a way of computing local arithmetic degrees in one common location namely the value group of $\cpt$--in particular local arithmetic degrees arising from different geometric data of $(\xan_E,\xan_K, \xan_K\to\xan_E)$ may be compared in one common location (this becomes important for formulating and proving \cite[Corollary 3.12]{mochizuki-iut3}).
\eenum
\erem
In many arithmetic contexts, one can take $F=\cpt$ and hence,
in particular one has the following corollary of Theorem~\ref{thm:arith-hol-strs}:

\bcor\label{cor:labeling-tempered}
Let $X/E$ be a geometrically connected, smooth quasi-projective variety over a $p$-adic field $E$ and $F=\cpt$.
Then the natural function $$(E\into K,K^\flat\isom F=\cpt)\mapsto \pitk{X/E}{K}$$ from the set of inequivalent untilts of $\C_p^\flat$ to the topological isomorphism class of $\pit{X/E}$
provides a distinguished  collection of distinctly labeled isomorphs  $$\left\{\pitk{X/E}{K}: (E\into K,K^\flat\isom F=\cpt) \in \syfe \right\}$$ of the tempered fundamental group  $\pit{X/E}$.
\ecor

The above consideration can be applied to \'etale fundamental groups of geometrically connected, smooth quasi-projective varieties as follows. Let $X/E$ be a geometrically connected, smooth, quasi-projective variety over a $p$-adic field $E$. Then one has a natural homomorphism (\cite[Proposition 4.4.1]{andre03}, \cite[Section 2.1.4]{andre-book}):
$$\pit{X/E}\to\pi_1(X/E),$$ 
which is injective if $\dim(X)=1$, and in any dimension
its image is dense and  $\pi_1(X/E)$ is the profinite completion  of the image of the above map
$$\widehat{\pit{X/E}}=\pi_1(X/E).$$ 

Let $(K,K^\flat\isom F)$ be an untilt of $F$. I define
$$\pi_1(X/E)_{(K,K^\flat\isom F)}=\widehat{\pitk{X/E}{K}},$$
and for simplicity write this as $\pi_1(X/E)_{K}$, and call $\pi_1(X/E)_{K}$ the \emph{untilt of the \'etale fundamental group} $\pi_1(X/E)$ corresponding to the untilt $(K,K^\flat\isom F)$ (of $F\isom K^\flat$). 
Thus one has the notion of untilts of the \'etale fundamental group $\pi_1(X/E)$.

\bcor
Let $X/E$ be a geometrically connected, smooth quasi-projective variety over a $p$-adic field $E$.
Then the natural function $$(K,K^\flat\isom F)\mapsto \pi_1(X/E)_{K}$$ from the set of inequivalent untilts of $F$ to the topological isomorphism class of $\pi_1(X/E)$
provides a distinguished collection of distinctly labeled copies  $$\left\{\pi_1(X/E)_{K}: y \in \syfe \text{ a closed classical point  with residue field } K_y \right\}$$ of the \'etale fundamental group  $\pi_1(X/E)$.
\ecor

I have used untilts of algebraically closed, perfectoid fields of characteristic $p>0$ as a set of distinguishing labels for the copies of fundamental groups produced here. There is another natural, but perhaps less useful indexing set:

\newcommand{\sK}{\mathcal{K}}

\bcor
Let $E$ be a $p$-adic field, $X/E$ a geometrically connected, smooth, quasi-projective variety over $E$. Consider the set of topological isomorphism classes of algebraically closed, complete valued fields $K\supset E$ (isometric inclusions):  $$\sK_E=\left\{K: E\subset K, K=\kbh\right\}.$$

Then there is a natural function $K\mapsto \pi_1(X/E)_K$ from $\sK_E$ to the topological isomorphism class of the profinite group $\pi_1(X/E)$ given by considering the tempered fundamental group associated to the datum $(X,E\into K)$.
\ecor

\brem\label{re:archimedean-case}
There is a further aspect of this result which should be pointed out. One should view elements $K\in\sK_E$ as providing a topological variation of ambient (additive) structure   $K\supset \bE$ while keeping internal field structure of $\bE$ unchanged. \emph{Such variations exist because, unlike the number field case,  algebraically closed fields such as $\C_p$, are quite far from being topologically rigid.} 
This is in complete contrast with the archimedean case, where by the well-known theorem of Ostrowski \cite[Chap. 6, \S 6, Th\'eor\`eme 2]{bourbaki-alg-com}, one knows that the only algebraically closed field complete with respect to an archimedean valuation is isometric to $\C$. To put Ostrowski's Theorem differently: \emph{Any two  algebraically closed, \emph{archimedean perfectoid fields} (i.e. fields which are algebraically closed and complete with respect to an archimedean valuation) are isometric (and also isometric with $\C$) and hence such fields are topologically rigid.}
\erem

\section{Arithmetic Teichmuller Spaces}\label{se:arith-hol-spaces}

Let me briefly sketch how the above ideas may be used to construct arithmetic Teichmuller spaces described in \cite{joshi-teich}. 

\begin{defn}
Let $E$ be a $p$-adic field and let $X/E$ be a geometrically connected, smooth, quasi-projective variety over $E$. Then the \textit{arithmetic Teichmuller space $\fjxe$ of $X/E$} is a category whose objects are

$$
\left\lbrace (Y/E',(E'\into K, K^\flat \isom F), *_K:\sM(K)\to \yan_{E'})\;:\;
\Vectorstack{
{F\text{  is some alg. closed perf. field }char(F)=p>0, }
{ (E'\into K, K^\flat \isom F)  \text{ an untilt of } F}
{Y/E' \text{ hyp. geom. con. smooth,} E' \text{ a } p\text{-adic field,}}
{Y/E' \text{ tempered anabelomorphic to } X/E, \text{ i.e. }}
{\pit{Y/E'}\isom{\pit{X/E}},}
{ \dim(Y)=\dim(X).}
}
\right\rbrace
$$
\end{defn}

Here \textit{tempered anabelomorphic} (see \cite[Section 2]{joshi-anabelomorphy}) means $Y/E'$ and $X/E$ have isomorphic tempered fundamental groups (this condition is independent of the choice of base-points and the choice of arithmetic holomorphic structures). [By standard results \cite{andre-book},  tempered anabelomorphisms provide  anabelomorphisms i.e. isomorphisms of \'etale fundamental groups of $Y/E'$ and $X/E$.] Morphisms between objects of $\fjxe$ are defined in the obvious way. 

\brem The following remarks regarding this definition of $\fjxe$ will be important.
\benumlab
\item The category $\fjxe$ is, of course, quite a big category, and in practice one must (and one will always) restrict $Y/E'$ in some way to avoid inane pathologies.
\item  For example, if $X/E$ is a $K(\pi,1)$ space, then one can require $Y/E'$ to be $K(\pi,1)$ space for every object $(Y/E',(E'\into K, K^\flat \isom F), *_K:\sM(K)\to \yan_{E'})$.
\item  If $X/E$ is an hyperbolic curve over a $p$-adic field, then one restricts $Y/E'$ to be an hyperbolic curve for every object $(Y/E',(E'\into K, K^\flat \isom F), *_K:\sM(K)\to \yan_{E'})$ etc. 
\item \textit{Notably if $X/E$ is an hyperbolic curve, I will restrict $Y/E'$ to be an hyperbolic curve as well.}
\item Note that the field $E'$ need not be isomorphic to $E$--the existence of an anabelomorphism $\pit{Y/E'}\isom \pit{X/E}$ forces an anabelomorphism $G_{E'}\isom G_E$. 
\eenum
\erem

Another way of thinking about objects of the arithmetic Teichmuller space $\fjxe$  is this:
$$\left\lbrace Y/E', \text{ arith. hol. strs. on }Y/E': 
\Vectorstack{
{E'\text{ a }p \text{-adic field,}}
{Y/E' \text{ tempered anabelomorphic to }X/E}
{\dim(Y)=\dim(X).}
}
\right\rbrace
$$

In \cite{joshi-teich-rosetta}, $Y/E'$ equipped with an arithmetic holomorphic structure is called a \textit{holomorphoid on} $Y/E'$.

With this definition of $\fjxe$, for every object $(Y/E',(E'\into K, K^\flat \isom F), *_K:\sM(K)\to \yan_{E'})\in\fjxe$, one has an isomorphism of topological groups $\pit{Y/E'}\isom \pit{X/E}$ and  an isomorphism of the corresponding \'etale fundamental groups and this also provides an isomorphism of the absolute Galois groups $G_{E'}\isom G_E$. 

The arithmetic Teichmuller space $\fjxe$  should be considered to be an anabelian analog of the variation of mixed Hodge structures (\cite{joshi-teich}). The arithmetic Teichmuller space comes equipped with many functors to several categories of interest. This is detailed in \cite{joshi-teich}.

\brem\label{re:arithmetic-space-point}
The following remark will be useful in understanding the complicated nature of the arithmetic Teichmuller space defined above. Let $E$ be a $p$-adic field and let $X=\Spec(E)$. Then it makes sense to consider the arithmetic Teichmuller space $\fjxe$. 

To understand $\fjxe$, let $F$ be a perfectoid field of characteristic $p>0$. Then an arbitrary object of $\fjxe$ consists of $(Y/E', (E'\into K,K^\flat\isom F), *_K:\sM(K)\to X^{an}_{E}=\sM(E))$ the condition $\dim(Y)=\dim(X)$ and connectedness and smoothness assumption forces $Y$ to be a finite, smooth, geometrically connected $E'$ algebra i.e. $Y=\spec(E')$ with $$G_{E'}=\pi_1^{temp}(Y/E')=\pi_1^{temp}(X/E)=G_E.$$ 

Thus for $X=\Spec(E)$, the objects of $\fjxe$ consists  of closed classical points of all Fargues-Fontaine curves $\syfep$ for all perfectoid fields $F$ of characteristic $p>0$ for all $p$-adic fields $E'$ that $$G_{E'}\isom G_E,$$
and notably for because of the canonical action of $G_E$ (resp. $G_{E'})$ on $\syfe$ (resp. $\syfep$) one can write this as
$$\syfep\curvearrowleft G_{E'} \isom G_E \curvearrowright\syfe.$$
Notably in terms of the complete Fargues-Fontaine curves one has an anabelomorphism of schemes (for the usual choice of basepoints on either sides)
$$\pi_1^{et}(\sxfep)=G_{E'}\isom G_E=\pi_1(\sxfe).$$
\textit{If $E',E$ are not isomorphic as $p$-adic fields, then by the main theorem of \cite{joshi-gconj}, this anabelomorphism need not arise from an isomorphism $\sxfep\isom\sxfe$ of schemes.}

Thus even for a simple case $X=\spec(E)$ of a point, the space $\fjxe$ is an highly non-trivial anabelian object which ``anabelomorphically glues'' together many Fargues-Fontaine curves  by means non-trivial anabelomorphisms between them.

One may also similarly describe the arithmetic Teichmuller space $\fJ(\P^n/E)$.
\erem

I will say that $\fjxe$ is a \textit{connected arithmetic Teichmuller space} if for every object $(Y/E',(E'\into K, K^\flat \isom F), *_K:\sM(K)\to \yan_{E'})\in\fjxe$, one has an isomorphism  $X\isom Y$ of schemes over $\Z$. Notably there is a tautological ``connected component,'' namely the full subcategory of $\fjxe$ whose class of objects is $$\left\lbrace (X/E,(E\into K, K^\flat \isom F), *_K:\sM(K)\to \xan_{E}) : F \text{ is some alg. closed perfectoid of } char(F)=p>0\right\rbrace.$$
The relationship between connectedness of arithmetic Teichmuller spaces and the absolute Grothendieck conjecture is summarized in the following:

\bthm\label{th:groth-and-connectedness} Let $E$ be a $p$-adic field and let $X/E$ be a geometrically connected, smooth, quasi-projective variety over $E$. Then 
\benumlab
\item $\fjxe$ is a connected arithmetic Teichmuller space if and only if the absolute Grothendieck conjecture holds true for any pair of objects of $\fjxe$. 
\item In particular if $X/E$ is additionally an hyperbolic curve of strict Belyi Type  \cite{mochizuki07} then the arithmetic Teichmuller space $\fjxe$ is connected.
\eenum
\ethm
\bp 
I will not recall the definition of hyperbolic curves of strict Belyi type and readers are referred to \cite{mochizuki07} for this notion and its properties.

The first assertion is clear from the definition of connectedness given above and the statement of the absolute Grothendieck conjecture. The second assertion is a consequence of the well-known result of \cite[Corollary 2.3]{mochizuki07} which shows that the absolute Grothendieck conjecture holds true for hyperbolic curves of strict Belyi type over $p$-adic fields, so the second assertion follows from the first assertion.
\ep

\ENDDOCUMENT

\section{Deformations of Number Fields}
The considerations presented in the previous sections deals mostly with $p$-adic fields for a fixed prime $p$. Let me now indicate how my ideas  can be applied to providing arithmetic deformations of number fields. This is an idea due to Shinichi Mochizuki, but I arrived at it quite differently and my presentation of this is based on the ideas presented in the previous sections, and the framework I present allows me to demonstrate the Mochizuki's claims with complete clarity.

Let me make a definition which will allow us to treat archimedean primes on par with non-archimedean primes. For treating archimedean valuations essentially on par with non-archimedean valuations. An \textit{archimedean algebraically closed perfectoid field} is a valued field $(K,\abs{-}_K)$ which is isometric with $(\C,\abs{-}_{\C})$. Let $(K,\abs{-}_K)$ be an algebraically closed, archimedean perfectoid field. Define $(K^\flat,\abs{-}_{K^\flat})=(K,\abs{}_K)$ in particular one has the  \textit{tilt of complex numbers} $\C$ given by  $(\C^\flat,\abs{-}_{\C^\flat})=(\C,\abs{-}_{\C})$. An untilt of $\C^\flat$ is a archimedean, algebraically closed perfectoid field $K$ and an isometry $K^\flat\isom \C^\flat$.

\newcommand{\V}{\mathbb{V}}
\newcommand{\lvbh}{\widehat{\overline{L}}_v}
I want to introduce the notion of arithmetic deformations of a number fields. Now let $L$ be a number field, $\V_L$ be the set of inequivalent non-trivial valuations of $L$. Let $v\in\V_L$ be a valuation of $L$. Let $L_v$ be the completion of $L$ at $v$, let $\lvbh$ be a fixed algebraic closure of $L_v$. Let $\lvbh^\flat$ be its tilt (because of the above definition, this makes sense for every $v\in\V_L$). For every valuation of $L$. An \textit{arithmeticoid of $L$} or an \textit{arithmetic deformation of $L$} is the following data. For each $v\in\V_L$ one is given an untilt $(L_v\into K_v,K_v^\flat\isom \lvbh^\flat)$. This local data is required to satisfy the following global arithmetic requirement (the validity of the product formula): for all $0\neq x\in L$, one has
\be 
 \prod_{v\in\V_L}\abs{x}_{K_v}=1.
\ee
I will say that two arithmeticoids of $L$ are  isomorphic if one has for every $v\in\V_L$ an isomorphism of untilts $(L_v\into K_v,K_v^\flat\isom \lvbh^\flat)\isom (L_v \into K_v',K_v^{'\flat}\isom \lvbh^\flat)$.

The following is a tautology:
\bthm 
Let $L$ be any number field. Then the class of arithmeticoids of $L$ is non-empty.
\ethm
\bp For each prime $v$ of $L$ choose an algebraic closure $\bL_v$ of $L_v$ and
let $\lvbh$ be the completion of $\bL_v$. Then one has a natural untilt $(L_v\into \widehat{\bL}_b,\lvbh^{\flat}=\lvbh^{\flat})$ and thus one has an an arithmeticoid 
\ep

It is clear that many  non-isomorphic arithmeticoids of $L$ exists. \textit{Let me explain why this a precise way of describing  Mochizuki's idea of dismantling additive and multiplicative structure of number fields.} For each non-archimedean prime $v$, Fargues-Fontaine curves $\syqp,\sxqp$\cite{fargues-fontaine} can be described additively (via a one dimensional Lubin-Tate formal group) and also multiplicatively using the multiplicative formal group $\hgm/\Z_p$. In the multiplicative view (compatible with Mochizuki's multiplicative monoid point of view in \iut), the additive structure of the perfectoid field $K_v$ varies while the multiplicative structure  $\invlim_{x\mapsto x^p}\hgm(\O_{K_v})\isom \hgm(\O_{\lvbh^\flat})$ remains fixed since one is given an isometry $K_v^\flat \isom \lvbh^\flat$. So inequivalent arithmeticoids of $L$ provide distinct views of $L$ by dismantling the ring structure of $L$ as claimed by Mochizuki. On the other hand each arithmeticoid of $L$ provides a distinct version of $L$ in which one has compatibility of multiplicative structures $L^*\to L_v^*$ but the additive structure are related only through the $p$-adic logarithm via the surjections $\hgm(\O_{\lvbh^\flat})\onto K_v$ (for each $v$.

\ENDDOCUMENT
\begin{thebibliography}{49}
\providecommand{\natexlab}[1]{#1}
\providecommand{\url}[1]{\texttt{#1}}
\expandafter\ifx\csname urlstyle\endcsname\relax
  \providecommand{\doi}[1]{doi: #1}\else
  \providecommand{\doi}{doi: \begingroup \urlstyle{rm}\Url}\fi

\bibitem[Andr\'{e}(2003)]{andre-book}
Yves Andr\'{e}.
\newblock \emph{Period mappings and differential equations. {F}rom {$\Bbb C$}
  to {$\Bbb C_p$}}, volume~12 of \emph{MSJ Memoirs}.
\newblock Mathematical Society of Japan, Tokyo, 2003.
\newblock T\^{o}hoku-Hokkaid\^{o} lectures in arithmetic geometry, With
  appendices by F. Kato and N. Tsuzuki.

\bibitem[Andr\'e(2003)]{andre03}
Yves Andr\'e.
\newblock On a geometric description of $\gal({\bQ}_p/{\Q}_p)$ and a $p$-adic
  avatar of $\widehat{GT}$.
\newblock \emph{Duke Math. Journal}, 119\penalty0 (1):\penalty0 1--39, 2003.

\bibitem[Beilinson and Schechtman(1988)]{beilinson88}
A.~A. Beilinson and V.~V. Schechtman.
\newblock Determinant bundles and {V}irasoro algebras.
\newblock \emph{Comm. Math. Phys.}, 118\penalty0 (4):\penalty0 651--701, 1988.
\newblock URL \url{http://projecteuclid.org/euclid.cmp/1104162170}.

\bibitem[Berkovich(1990)]{berkovich-book}
Vladimir~G. Berkovich.
\newblock \emph{Spectral theory and analytic geometry over non-{A}rchimedean
  fields}, volume~33 of \emph{Mathematical Surveys and Monographs}.
\newblock American Mathematical Society, Providence, RI, 1990.

\bibitem[Bourbaki(1985)]{bourbaki-alg-com}
Nicolas Bourbaki.
\newblock \emph{\'{E}l\'{e}ments de math\'{e}matique}.
\newblock Masson, Paris, 1985.
\newblock ISBN 2-225-80269-6.
\newblock Alg\`ebre commutative. Chapitres 5 \`a 7. [Commutative algebra.
  Chapters 5--7], Reprint.

\bibitem[Fargues and Fontaine(2018)]{fargues-fontaine}
Laurent Fargues and Jean-Marc Fontaine.
\newblock Courbes et fibr\'{e}s vectoriels en th\'{e}orie de {H}odge
  {$p$}-adique.
\newblock \emph{Ast\'{e}risque}, \penalty0 (406):\penalty0 xiii+382, 2018.
\newblock ISSN 0303-1179.
\newblock With a preface by Pierre Colmez.

\bibitem[Frenkel and Ben-Zvi(2001)]{frenkel01-book}
Edward Frenkel and David Ben-Zvi.
\newblock \emph{Vertex algebras and algebraic curves}, volume~88 of
  \emph{Mathematical Surveys and Monographs}.
\newblock American Mathematical Society, Providence, RI, 2001.
\newblock URL \url{https://doi.org/10.1090/surv/088}.

\bibitem[Grothendieck(1971)]{grothendieck1971a}
A.~Grothendieck.
\newblock \emph{Rev$\hat e$tment \'etale et groupe fondamental, S\'eminaire de
  G\'eom\'etrie alg\'ebrique du Bois-Marie 1960-61}.
\newblock Number 224 in Springer Lecture Notes in Math. Springer-Verlag, 1971.

\bibitem[{John Keats}()]{keats}
{John Keats}.
\newblock \emph{The poetical works of John Keats}.
\newblock Frederick Warne \& Co.

\bibitem[Joshi(2019)]{joshi-formal-groups}
Kirti Joshi.
\newblock Mochizuki's anabelian variation of ring structures and formal groups.
\newblock 2019.
\newblock URL \url{https://arxiv.org/pdf/1906.06840.pdf}.
\newblock https://arxiv.org/abs/1906.06840.

\bibitem[Joshi(2020{\natexlab{a}})]{joshi-anabelomorphy}
Kirti Joshi.
\newblock On {M}ochizuki's idea of anabelomorphy and applications.
\newblock 2020{\natexlab{a}}.
\newblock URL \url{https://arxiv.org/pdf/2003.01890.pdf}.

\bibitem[Joshi(2020{\natexlab{b}})]{joshi-gconj}
Kirti Joshi.
\newblock The absolute {G}rothendieck conjecture is false for
  {F}argues-{F}ontaine curves.
\newblock 2020{\natexlab{b}}.
\newblock URL \url{https://arxiv.org/pdf/2008.01228.pdf}.
\newblock Preprint.

\bibitem[Joshi(2020{\natexlab{c}})]{joshi-untilts}
Kirti Joshi.
\newblock Untilts of fundamental groups: construction of labeled isomorphs of
  fundamental groups.
\newblock 2020{\natexlab{c}}.
\newblock URL \url{https://arxiv.org/pdf/2010.05748.pdf}.

\bibitem[Joshi(2021{\natexlab{a}})]{joshi-teich}
Kirti Joshi.
\newblock Construction of {A}rithmetic {T}eichmuller {S}paces and applications.
\newblock 2021{\natexlab{a}}.
\newblock URL \url{https://arxiv.org/abs/2106.11452}.

\bibitem[Joshi(2021{\natexlab{b}})]{joshi-teich-estimates}
Kirti Joshi.
\newblock Construction of {A}rithmetic {T}eichmuller {S}paces and some
  applications {II}: {T}owards {D}iophantine estimates.
\newblock 2021{\natexlab{b}}.
\newblock URL \url{https://arxiv.org/pdf/2111.04890.pdf}.

\bibitem[Joshi(2021{\natexlab{c}})]{joshi-teich-summary-comments}
Kirti Joshi.
\newblock Comments on {A}rithmetic {T}eichmuller {T}heory.
\newblock 2021{\natexlab{c}}.
\newblock URL \url{https://arxiv.org/pdf/2111.06771.pdf}.

\bibitem[Joshi(2022)]{joshi-teich-rosetta}
Kirti Joshi.
\newblock Construction of {A}rithmetic {T}eichmuller {S}paces and applications
  {III}.
\newblock \emph{In preparation}, 2022.

\bibitem[Kedlaya and Temkin(2018)]{kedlaya18}
Kiran~S. Kedlaya and Michael Temkin.
\newblock Endomorphisms of power series fields and residue fields of
  {F}argues-{F}ontaine curves.
\newblock \emph{Proc. Amer. Math. Soc.}, 146\penalty0 (2):\penalty0 489--495,
  2018.
\newblock URL \url{https://doi.org/10.1090/proc/13818}.

\bibitem[Kontsevich(1987)]{kontsevich87}
Maxim Kontsevich.
\newblock Virasoro algebra and {T}eichmuller spaces.
\newblock \emph{Functional analysis and its applications}, 21\penalty0
  (2):\penalty0 156, 1987.
\newblock \doi{10.1007/BF01078034}.

\bibitem[Lehto(1987)]{lehto-book}
Olli Lehto.
\newblock \emph{Univalent functions and {T}eichm\"{u}ller spaces}, volume 109
  of \emph{Graduate Texts in Mathematics}.
\newblock Springer-Verlag, New York, 1987.
\newblock URL \url{https://doi.org/10.1007/978-1-4613-8652-0}.

\bibitem[Lepage(2010)]{lepage-thesis}
Emmanuel Lepage.
\newblock Géométrie anabélienne tempérée ({PhD}. {T}hesis), 2010.
\newblock URL \url{https://arxiv.org/abs/1004.2150}.

\bibitem[Matignon and Reversat(1984)]{matignon84}
Michel Matignon and Marc Reversat.
\newblock Sous-corps ferm\'{e}s d'un corps valu\'{e}.
\newblock \emph{J. Algebra}, 90\penalty0 (2):\penalty0 491--515, 1984.
\newblock URL \url{https://doi.org/10.1016/0021-8693(84)90186-8}.

\bibitem[Mirzakhani(2007)]{mirzakhani07}
Maryam Mirzakhani.
\newblock Simple geodesics and weil-petersson volumes of moduli spaces of
  bordered {R}iemann surfaces.
\newblock \emph{Invent. {M}ath.}, 167:\penalty0 179--222, 2007.
\newblock \doi{https://doi.org/10.1007/s00222-006-0013-2}.

\bibitem[Mochizuki()]{mochizuki-ss-rejoinder}
Shinichi Mochizuki.
\newblock Report on discussions, held during the period {M}arch 15--20, 2018,
  concerning {I}nter-{U}niversal {T}eichmuller {T}heory.
\newblock URL \url{https://www.kurims.kyoto-u.ac.jp/~motizuki/Rpt2018.pdf}.

\bibitem[Mochizuki(1996{\natexlab{a}})]{mochizuki-ordinary}
Shinichi Mochizuki.
\newblock A theory of ordinary {$p$}-adic curves.
\newblock \emph{Publ. Res. Inst. Math. Sci.}, 32\penalty0 (6):\penalty0
  957--1152, 1996{\natexlab{a}}.
\newblock \doi{10.2977/prims/1195145686}.

\bibitem[Mochizuki(1996{\natexlab{b}})]{mochizuki96-gconj}
Shinichi Mochizuki.
\newblock The profinite {G}rothendieck conjecture for closed hyperbolic curves
  over number fields.
\newblock \emph{J. Math. Sci. Univ. Tokyo}, 3\penalty0 (3):\penalty0 571--627,
  1996{\natexlab{b}}.

\bibitem[Mochizuki(1999{\natexlab{a}})]{mochizuki-HA}
Shinichi Mochizuki.
\newblock The {H}odge-{A}rakelov theory of elliptic curves: Global
  discretization of local hodge theories.
\newblock 1999{\natexlab{a}}.
\newblock URL
  \url{https://www.kurims.kyoto-u.ac.jp/~motizuki/The%20Hodge-Arakelov%20Theory%20of%20Elliptic%20Curves.pdf}.

\bibitem[Mochizuki(1999{\natexlab{b}})]{mochizuki-teichmuller-book}
Shinichi Mochizuki.
\newblock \emph{Foundations of {$p$}-adic {T}eichm\"{u}ller theory}, volume~11
  of \emph{AMS/IP Studies in Advanced Mathematics}.
\newblock American Mathematical Society, Providence, RI; International Press,
  Cambridge, MA, 1999{\natexlab{b}}.

\bibitem[Mochizuki(2007)]{mochizuki07}
Shinichi Mochizuki.
\newblock A combinatorial version of the {G}rothendieck conjecture.
\newblock \emph{Tohoku Math. J. (2)}, 59\penalty0 (3):\penalty0 455--479, 2007.
\newblock ISSN 0040-8735.
\newblock URL \url{http://projecteuclid.org/euclid.tmj/1192117988}.

\bibitem[Mochizuki(2008)]{mochizuki-frobenioid1}
Shinichi Mochizuki.
\newblock The geometry of {F}robenioids. {I}. {T}he general theory.
\newblock \emph{Kyushu J. Math.}, 62\penalty0 (2):\penalty0 293--400, 2008.
\newblock \doi{10.2206/kyushujm.62.293}.
\newblock URL \url{https://doi.org/10.2206/kyushujm.62.293}.

\bibitem[Mochizuki(2012)]{mochizuki-topics1}
Shinichi Mochizuki.
\newblock Topics in absolute anabelian geometry {I}: generalities.
\newblock \emph{J. Math. Sci. Univ. Tokyo}, 19\penalty0 (2):\penalty0 139--242,
  2012.

\bibitem[Mochizuki(2013)]{mochizuki-topics2}
Shinichi Mochizuki.
\newblock Topics in absolute anabelian geometry {II}: decomposition groups and
  endomorphisms.
\newblock \emph{J. Math. Sci. Univ. Tokyo}, 20\penalty0 (2):\penalty0 171--269,
  2013.

\bibitem[Mochizuki(2015)]{mochizuki-topics3}
Shinichi Mochizuki.
\newblock Topics in absolute anabelian geometry {III}: global reconstruction
  algorithms.
\newblock \emph{J. Math. Sci. Univ. Tokyo}, 22\penalty0 (4):\penalty0
  939--1156, 2015.
\newblock ISSN 1340-5705.

\bibitem[Mochizuki(2021{\natexlab{a}})]{mochizuki-iut1}
Shinichi Mochizuki.
\newblock Inter-{U}niversal {T}eichmuller theory {I}: construction of {H}odge
  {T}heaters.
\newblock \emph{Publ. Res. Inst. Math. Sci.}, 57\penalty0 (1/2):\penalty0
  3--207, 2021{\natexlab{a}}.
\newblock URL
  \url{http://www.kurims.kyoto-u.ac.jp/~motizuki/Inter-universal%20Teichmuller%20Theory%20I.pdf}.

\bibitem[Mochizuki(2021{\natexlab{b}})]{mochizuki-iut2}
Shinichi Mochizuki.
\newblock Inter-{U}niversal {T}eichmuller {T}heory {II}: {H}odge-{A}rakelov
  {T}heoretic {E}valuations.
\newblock \emph{Publ. Res. Inst. Math. Sci.}, 57\penalty0 (1/2):\penalty0
  209--401, 2021{\natexlab{b}}.
\newblock URL
  \url{http://www.kurims.kyoto-u.ac.jp/~motizuki/Inter-universal%20Teichmuller%20Theory%20II.pdf}.

\bibitem[Mochizuki(2021{\natexlab{c}})]{mochizuki-iut3}
Shinichi Mochizuki.
\newblock Inter-{U}niversal {T}eichmuller {T}heory {III}: canonical splittings
  of the log-theta lattice.
\newblock \emph{Publ. Res. Inst. Math. Sci.}, 57\penalty0 (1/2):\penalty0
  403--626, 2021{\natexlab{c}}.
\newblock URL
  \url{http://www.kurims.kyoto-u.ac.jp/~motizuki/Inter-universal%20Teichmuller%20Theory%20III.pdf}.

\bibitem[Mochizuki(2021{\natexlab{d}})]{mochizuki-iut4}
Shinichi Mochizuki.
\newblock Inter-{U}niversal {T}eichmuller {T}heory {IV}: Log-volume
  computations and set theoretic foundations.
\newblock \emph{Publ. Res. Inst. Math. Sci.}, 57\penalty0 (1/2):\penalty0
  627--723, 2021{\natexlab{d}}.
\newblock URL
  \url{http://www.kurims.kyoto-u.ac.jp/~motizuki/Inter-universal%20Teichmuller%20Theory%20IV.pdf}.

\bibitem[Mochizuki(2022)]{mochizuki-essential-logic}
Shinichi Mochizuki.
\newblock On the essential logical structure of {I}nter-{U}niversal
  {T}eichmuller {T}heory in terms of logical and “$\wedge$”/logical or
  “$\vee$” relations: Report on the occasion of the publication of the four
  main papers on {I}nter-{U}niversal {T}eichmuller {T}heory.
\newblock 2022.
\newblock URL
  \url{https://www.kurims.kyoto-u.ac.jp/~motizuki/Essential%20Logical%20Structure%20of%20Inter-universal%20{T}eichmuller%20Theory.pdf}.

\bibitem[Morrow(2019)]{morrow-bourbaki}
Matthew Morrow.
\newblock The {F}argues-{F}ontaine curve and diamonds [d'apr\`es {F}argues,
  {F}ontaine, and {S}cholze].
\newblock \emph{Ast\'{e}risque}, \penalty0 (414, S\'{e}minaire Bourbaki. Vol.
  2017/2018. Expos\'{e}s 1136--1150):\penalty0 Exp. No. 1150, 533--572, 2019.
\newblock URL \url{https://doi.org/10.24033/ast.1094}.

\bibitem[Mumford(1983)]{mumford-theta1}
David Mumford.
\newblock \emph{Tata Lectures on Theta I}.
\newblock Birkhäuser Boston, Boston, MA, 1983.

\bibitem[Scholze(2012)]{scholze12-perfectoid-ihes}
Peter Scholze.
\newblock Perfectoid spaces.
\newblock \emph{Publ. Math. Inst. Hautes \'{E}tudes Sci.}, 116:\penalty0
  245--313, 2012.
\newblock URL \url{https://doi.org/10.1007/s10240-012-0042-x}.

\bibitem[Scholze(2017)]{scholze-diamonds}
Peter Scholze.
\newblock \emph{\'Etale cohomology of Diamonds}.
\newblock 2017.
\newblock URL \url{https://arxiv.org/pdf/1709.07343}.

\bibitem[Scholze(2021)]{scholze-review}
Peter Scholze.
\newblock Review of {M}ochizuki's paper: {I}nter-{U}niversal {T}eichmüller
  {T}heory. {I},{II},{III},{IV}.
\newblock \emph{zbMath Open (formerly Zentralblatt Math)}, 2021.
\newblock URL \url{https://zbmath.org/pdf/07317908.pdf}.

\bibitem[Scholze and Stix()]{scholze-stix}
Peter Scholze and Jakob Stix.
\newblock Why $abc$ is still a conjecture.
\newblock URL \url{http://www.kurims.kyoto-u.ac.jp/~motizuki/SS2018-08.pdf}.

\bibitem[Scholze and Weinstein()]{scholze-weinstein-book}
Peter Scholze and Jared Weinstein.
\newblock \emph{Berkeley lectures on $p$-adic geometry}.
\newblock URL \url{http://www.math.uni-bonn.de/people/scholze/Berkeley.pdf}.

\bibitem[Scholze and Weinstein(2020)]{scholze-weinstein}
Peter Scholze and Jared Weinstein.
\newblock \emph{Berkeley lectures on {$p$}-adic geometry}, volume 207 of
  \emph{Annals of Mathematics Studies}.
\newblock Princeton University Press, Princeton, NJ, 2020.

\bibitem[Tamagawa(1997)]{tamagawa97-gconj}
Akio Tamagawa.
\newblock The {G}rothendieck conjecture for affine curves.
\newblock \emph{Compositio Math.}, 109\penalty0 (2):\penalty0 135--194, 1997.
\newblock URL \url{https://doi.org/10.1023/A:1000114400142}.

\bibitem[Tan(2018)]{fucheng}
Fucheng Tan.
\newblock Note on {I}nter-{U}niversal {T}eichmuller {T}heory.
\newblock 2018.

\bibitem[Wright(2019)]{wright19}
Alex Wright.
\newblock A tour through {M}irzakhani's work on moduli spaces of {R}iemann
  surfaces, 2019.
\newblock URL \url{https://arxiv.org/abs/1905.01753}.

\end{thebibliography}
